\documentclass[
	english
]{scrartcl}
\usepackage[T1]{fontenc}
\usepackage[utf8]{inputenc}
\usepackage[english]{babel}
\usepackage[inline,shortlabels]{enumitem}
\usepackage{amsmath}
\usepackage{amssymb}
\usepackage{amsthm}
\usepackage[noadjust]{cite}
\usepackage[a4paper,total={6in,9in}]{geometry}
\usepackage[colorlinks,citecolor=blue]{hyperref}
\usepackage[capitalize,noabbrev,nameinlink]{cleveref}
\usepackage{xcolor}
\usepackage{longtable}
\usepackage{tikz,pgfplots,pgfplotstable}
\usetikzlibrary{patterns}
\usepackage{algorithmic}
\usepackage{marginnote}
\usepackage[labelformat=simple]{subcaption}

\setlist{font=\normalfont,topsep=1ex,parsep=0ex}

\setlist[enumerate]{label=(\alph*)}

\usepackage{changes}

\numberwithin{equation}{section}
\numberwithin{table}{section}
\numberwithin{figure}{section}

\usepackage[figure]{hypcap}
\usepackage{graphicx}
\graphicspath{{./img/}}

\usepackage{preamble}
\usepackage{pdflscape}
\usepackage{booktabs}

\crefname{figure}{Figure}{Figures}
\crefname{table}{Table}{Tables}
\crefname{assumption}{Assumption}{Assumptions}
\Crefname{ALC@unique}{Step}{Steps}

\newlist{alglist}{enumerate}{1}
\setlist[alglist]{topsep=1ex,parsep=0ex,leftmargin=*,label=\textbf{Step~\arabic*.}}

\let\eqref\labelcref

\allowdisplaybreaks

\newcommand\norm[1]{\left\Vert#1\right\Vert}
\newcommand\nnorm[1]{\Vert#1\Vert}

\newcommand\N{\mathbb{N}}
\newcommand\R{\mathbb{R}}

\newcommand\calF{\mathcal F}
\newcommand\calG{\mathcal G}
\newcommand\calH{\mathcal H}

\newcommand\tto{\rightrightarrows}

\newcommand{\diag}{\operatorname{diag}}

\newcommand{\dom}{\operatorname{dom}}

\newcommand{\gph}{\operatorname{gph}}

\DeclareMathOperator*{\argmin}{\operatorname{argmin}}

\DeclareMathAlphabet{\mathpzc}{OT1}{pzc}{m}{it}
\newcommand\oo{\mathpzc{o}}
\newcommand\OO{\mathpzc{O}}

\newtheorem{theorem}{Theorem}[section]
\newtheorem{lemma}[theorem]{Lemma}

\newtheorem{corollary}[theorem]{Corollary}
\newtheorem{remark}[theorem]{Remark}
\newtheorem{definition}[theorem]{Definition}
\newtheorem{example}[theorem]{Example}

\newtheorem{algorithm}[theorem]{Algorithm}

\definecolor{mygreen}{rgb}{0.0,0.7,0.0}
\definecolor{mybrown}{rgb}{0.5,0.5,0.0}

\begin{document}

\title{A fresh look at nonsmooth Levenberg--Marquardt methods with applications to bilevel optimization}
%\subtitle{}
\author{%
	Lateef O.\ Jolaoso%
	\footnote{%
		University of Southampton,
		School of Mathematics,
		Southampton SO17 1BJ,
		United Kingdom,
		\email{l.o.jolaoso@soton.ac.uk},
		\orcid{0000-0002-4838-7465}
		}%
	\and
	Patrick Mehlitz%
	\footnote{%
		University of Duisburg-Essen,
		Faculty of Mathematics,
		45127 Essen,
		Germany,
		\email{patrick.mehlitz@uni-due.de},
		\orcid{0000-0002-9355-850X}%
		}%
	\and
	Alain B.\ Zemkoho%
	\footnote{%
		University of Southampton,
		School of Mathematics,
		Southampton SO17 1BJ,
		United Kingdom,
		\email{a.b.zemkoho@soton.ac.uk},
		\orcid{0000-0003-1265-4178}
		}%
	}

% \date{\today}
\publishers{}
\maketitle

\begin{abstract}
	In this paper, we revisit the classical problem of solving over-determined
	systems of nonsmooth equations numerically. 
	We suggest a nonsmooth Levenberg--Marquardt method for its solution which,
	in contrast to the existing literature, does not require local Lipschitzness
	of the data functions. This is possible when using Newton-differentiability
	instead of semismoothness as the underlying tool of generalized differentiation.
	Conditions for local fast convergence of the method are given.
	Afterwards, in the context of over-determined mixed nonlinear complementarity systems,
	our findings are applied, and globalized solution methods, based on a residual
	induced by the maximum and the Fischer--Burmeister function, respectively, are constructed.
	The assumptions for local fast convergence are worked out and compared.
	Finally, these methods are applied for the numerical solution of bilevel optimization
	problems. We recall the derivation of a stationarity condition taking the shape of
	an over-determined mixed nonlinear complementarity system involving a penalty parameter, 
	formulate assumptions for local fast convergence of our solution methods explicitly, 
	and present results of numerical experiments. 
	Particularly, we investigate whether the treatment of the appearing penalty parameter
	as an additional variable is beneficial or not.
\end{abstract}


\begin{msc}	
	\mscLink{49M05}, 
	\mscLink{49M15},
	\mscLink{90C30},
	\mscLink{90C33}
\end{msc}

\section{Introduction}\label{sec:introduction}

Nowadays, the numerical solution of nonsmooth equations is a classical topic in
computational mathematics. Journeying thirty years back to the past, Qi and colleagues
developed the celebrated nonsmooth Newton method, see \cite{Qi1993,QiSun1993}, based
on the notion of semismoothness for locally Lipschitzian functions
which is due to Mifflin, see \cite{Mifflin1977}.
Keeping the convincing local convergence properties of Newton-type methods in mind,
this allowed for a fast numerical solution of Karush--Kuhn--Tucker systems,
associated with constrained nonlinear optimization problems with inequality constraints, without the
need of smoothing or relaxing the involved complementarity slackness condition.
Indeed, the latter can be reformulated as a nonsmooth equation with the aid of so
called NCP-functions,
see \cite{Galantai2012,KanzowYamashitaFukushima1997,SunQi1999} for an overview,
where NCP abbreviates nonlinear complementarity problem.
Recall that a continuous function $\varphi\colon\R^2\to\R$ is referred to as an
NCP-function whenever 
\[
	\forall (a,b)\in\R^2\colon\quad
	\varphi(a,b)=0
	\quad\Longleftrightarrow
	a\leq 0,\,b\leq 0,\,ab=0
\]
is valid.
A rather prominent example of an NCP-function is the famous Fischer--Burmeister (FB) function
$\varphi_\textup{FB}\colon\R^2\to\R$ given by
\begin{equation}\label{eq:FB_function}
	\forall (a,b)\in\R^2\colon\quad
	\varphi_\textup{FB}(a,b):=a+b+\sqrt{a^2+b^2},
\end{equation}
see \cite{Fischer1992} for its origin. Another popular NCP-function is
the maximum function $(a,b)\mapsto\max(a,b)$. 
We refer the reader to 
\cite{DeLucaFacchineiKanzow1996,DeLucaFacchineiKanzow2000,FacchineiKanzow1997,Fischer1992}
for further early references dealing with semismooth Newton-type methods to solve nonlinear
complementarity systems.
Let us note that even prior to the development of semismooth Newton methods thirty years ago, 
the numerical solution of nonsmooth equations via Newton-type algorithms
has been addressed in the literature, see e.g.
\cite{KojimaShindo1986,Kummer1988,Kummer1992,Pang1990,Robinson1994}.
We refer the interested reader to the monographs \cite{FacchineiPang2003,IzmailovSolodov2014,KlatteKummer2002}
for a convincing overview.

Let us mention some popular directions of research which have been enriched or motivated by the theory
of semismooth Newton methods. 
Exemplary, this theory has been extended to Levenberg--Marquardt (LM) methods in order to allow for the treatment
of over-determined or irregular systems of equations, 
see e.g.\ \cite{FacchineiKanzow1997,Jiang1999,KanzowPetra2004,MaWang2009,QiWuZhou2003}.
It also has been applied in order to solve so-called mixed nonlinear complementarity problems which are a
natural extension of nonlinear complementarity problems and may also involve pure (smooth) equations,
see e.g.\ \cite{ChenMangasarian1996,DarynkaIzmailovSolodov2005,HuangWang2012,Kanzow2000,LiFukushima2000,MonteiroPang1996}.
Recently, the variational notion of semismoothness* for set-valued mappings has been introduced in 
\cite{GfrererOutrata2021} in order to discuss the convergence properties of Newton's method for the numerical 
solution of so-called generalized equations given via set-valued operators.
There also exist infinite-dimensional extensions of nonsmooth Newton methods, 
see e.g.\ \cite{BrokateUlbrich2022,ChenNashedQi2000,HintermuellerItoKunisch2002,ItoKunisch2008,Schiela2008,Ulbrich2002,Ulbrich2011}.
In contrast to the finite-dimensional case, in infinite-dimensions, the underlying concept of semismoothness is
often replaced by so-called Newton-differentiability. 
In principle, this notion of generalized differentiation encapsulates all necessities to address the convergence analysis
of associated Newton-type methods properly. Noteworthy, it does not rely on the Lipschitzness of the underlying mapping.
In the recent paper \cite{HarderMehlitzWachsmuth2021}, this favorable advantage of Newton-derivatives has been used
to construct a Newton-type method for the numerical solution of a certain stationarity system associated with so-called
mathematical programs with complementarity constraints, as these stationarity conditions can be reformulated as a system
of discontinuous equations.

Let us point the reader's attention to the fact that some stationarity conditions for bilevel optimization problems,
see \cite{Dempe2002,DempeKalashnikovPerezValdesKalashnykova2015} for an introduction, can be reformulated as
(over-determined) systems of nonsmooth equations involving an (unknown) penalization parameter $\lambda$. 
This observation has been used in the recent papers 
\cite{FischerZemkohoZhou2021,FliegeTinZemkoho2021,TinZemkoho2021} in order to solve bilevel optimization problems
numerically. In \cite{FischerZemkohoZhou2021}, the authors used additional dummy variables in order to transform the
naturally over-determined system into a square system, and applied a globalized semismooth Newton method for its
computational solution. The papers \cite{FliegeTinZemkoho2021,TinZemkoho2021} directly tackle the over-determined
nonsmooth system with Gauss--Newton or LM methods. However, in order to deal with the nonsmoothness,
they either assume that strict complementarity holds or they smooth the complementarity condition.
In all three papers \cite{FischerZemkohoZhou2021,FliegeTinZemkoho2021,TinZemkoho2021}, it has been pointed out
that the choice of the penalty parameter $\lambda$ in the stationarity system is, on the one hand,
crucial for the success of the approach but, on the other hand, difficult to realize in practice.

The contributions in this paper touch several different aspects.
Motivated by the results from \cite{HarderMehlitzWachsmuth2021}, we study the numerical solution of over-determined
systems of nonsmooth equations with the aid of LM methods based on the notion of Newton-differentiability.
We show local superlinear convergence of this method under reasonable assumptions,
see \cref{thm:local_convergence_LM}.
Furthermore, we point out that even for local quadratic convergence, local Lipschitzness of the underlying mapping is not
necessary. To be more precise, a one-sided Lipschitz estimate, which is referred to as calmness in the literature, is
enough for that purpose. Thus, our theory applies to situations where the underlying mapping can be even discontinuous.
Our next step is the application of this method for the numerical solution of over-determined 
mixed nonlinear complementarity systems.
Here, we strike a classical path and reformulate the appearing complementarity conditions with the FB function
and the maximum function in order to obtain a system of nonlinear equations.
To globalize the associated Newton methods, we exploit the well-known fact that the squared norm of the residual
associated with the FB function is continuously differentiable, and make use of gradient steps
with respect to this squared norm in combination with an Armijo step size selection if LM
directions do not yield sufficient improvements.
We work out our abstract conditions guaranteeing local fast convergence for both approaches and state global convergence
results in \cref{thm:global_convergence_LM,thm:global_convergence_LM_FB}.
Furthermore, we compare these findings with the ones in the classical paper \cite{DeLucaFacchineiKanzow2000} where a related comparison
has been done in the context of semismooth Newton-type methods for square systems of nonsmooth equations.
Finally, we apply both methods in order to solve bilevel optimization problems via their over-determined stationarity
system of nonsmooth equations mentioned earlier. In contrast to \cite{FliegeTinZemkoho2021,TinZemkoho2021}, we neither
assume strict complementarity to hold nor do we smooth or relax the complementarity condition.
Some assumptions for local fast convergence are worked out, see \cref{thm:local_fast_convergence_obpp,thm:local_fast_convergence_obpp_FB}.
As extensive numerical testing of related approaches has been carried out in 
\cite{FischerZemkohoZhou2021,FliegeTinZemkoho2021,TinZemkoho2021}, our experiments in \cref{sec:computational_experiments}
focus on specific features of the developed algorithms.
Particularly, we comment on the possibility to treat the aforementioned penalty parameter $\lambda$ as a variable, 
investigate this approach numerically, and compare the outcome with the situation where $\lambda$ is handled as a parameter.

The remainder of the paper is organized as follows.
\cref{sec:preliminaries} provides an overview of the notation and preliminaries used in the manuscript.
Furthermore, we recall the notion of Newton-differentiability 
and present some associated calculus rules in \cref{sec:Newton_differentiability}.
A local nonsmooth LM method for Newton-differentiable mappings is discussed in \cref{sec:LM_local}.
\Cref{sec:LM} is dedicated to the derivation of globalized nonsmooth LM methods for the
numerical solution of over-determined mixed nonlinear complementarity systems.
First, we provide the analysis for the situation where the complementarity condition is reformulated with the
maximum function in \cref{sec:global_LM_max}.
Afterwards, in \cref{sec:global_LM_FB}, we comment on the differences which pop up when the FB function
is used instead.
The obtained theory is used in \cref{sec:optimistic_bpp} in order to solve bilevel optimization problems.
In \cref{sec:obpp_modeling}, we first recall the associated stationarity system of interest and characterize some
scenarios where it provides necessary optimality conditions. Second, we present different ways on how to model
these stationarity conditions as an over-determined mixed nonlinear complementarity system.
The numerical solution of these systems with the aid of nonsmooth LM methods is then investigated
in \cref{sec:computational_experiments} where results of computational experiments are evaluated.
The paper closes with some concluding remarks in \cref{sec:conclusions}.

\section{Notation and preliminaries}\label{sec:preliminaries}

\subsection{Notation}\label{sec:notation}

By $\N$, we denote the positive integers.
The set of all real matrices with $m\in\N$ rows and $n\in\N$ columns will be represented by $\R^{m\times n}$,
and $\mathbb O$ denotes the all-zero matrix of appropriate dimensions
while we use $\mathbb I^n$ for the identity matrix in $\R^{n\times n}$.
If $M\in\R^{m\times n}$ is a matrix and $I\subset\{1,\ldots,m\}$ is arbitrary,
then $M_I$ shall be the matrix which results from $M$ by deleting those
rows whose associated index does not belong to $I$.
Whenever the quadratic matrix $N\in\R^{n\times n}$ is symmetric,
$\lambda_\textup{min}(N)$ denotes the smallest eigenvalue of $N$.
For any $x\in\R^n$, $\diag(x)\in\R^{n\times n}$ is the diagonal matrix whose main diagonal is
given by $x$ and whose other entries vanish.

For arbitrary $n\in\N$, the space $\R^n$ will be equipped by the standard Euclidean inner product
as well as the Euclidean norm $\norm{\cdot}\colon\R^n\to\R$.
Further, we equip $\R^{m\times n}$ with the matrix norm induced by the Euclidean norm, i.e.,
with the spectral norm, and denote it by $\norm{\cdot}\colon\R^{m\times n}\to\R$ as well
as this cannot cause any confusion.
For arbitrary $x\in\R^n$ and $\varepsilon>0$, $\mathbb B_\varepsilon(x):=\{y\in\R^n\,|\,\norm{y-x}\leq\varepsilon\}$
represents the closed $\varepsilon$-ball around $x$.
Recall that a sequence $\{x^k\}_{k\in\N}\subset\R^n$ is said to converge superlinearly to some $\bar x\in\R^n$
whenever $\nnorm{x^{k+1}-\bar x}\in\oo(\nnorm{x^k-\bar x})$.
The convergence $x^k\to\bar x$ is said to be quadratic if $\nnorm{x^{k+1}-\bar x}\in\OO(\nnorm{x^k-\bar x}^2)$.

Frequently, for brevity of notation, we interpret tuples of vectors as a single block column vector, i.e.,
we exploit the identity
\[
	\forall i\in \{1,\ldots,\ell\},\,\forall x^i\in\R^{n_i}\colon\quad
	\R^{n_1}\times\ldots\times\R^{n_\ell}\ni(x^1,\ldots,x^\ell)
	\cong
	\begin{bmatrix}x^1\\\vdots\\x^\ell\end{bmatrix}\in\R^{n_1+\ldots+n_\ell}
\]
for any $\ell\in\N$ with $\ell\geq 2$ and $n_1,\ldots,n_\ell\in\N$.

A function $H\colon\R^n\to\R^m$ is called calm at $\bar x\in\R^n$ whenever there
are constants $C>0$ and $\varepsilon>0$ such that
\[
	\forall x\in\mathbb B_\varepsilon(\bar x)\colon\quad
	\nnorm{H(x)-H(\bar x)}\leq C\nnorm{x-\bar x}.
\]
We note that calmness can be seen as a one-sided local Lipschitz property which is,
in general, weaker than local Lipschitzness.
Indeed, one can easily check that there exist functions which are calm and
discontinuous at some fixed reference point.

Next, let $H\colon\R^n\to\R^m$ be differentiable at $\bar x\in\R^n$.
Then $H'(\bar x)\in\R^{m\times n}$ denotes the Jacobian of $H$ at $\bar x$.
Similarly, for a differentiable scalar-valued function $h\colon\R^n\to\R$,
$\nabla h(\bar x)\in\R^n$ is the gradient of $h$ at $\bar x$.
Clearly, $h'(\bar x)=\nabla h(\bar x)^\top$ by construction.
In the case where $h$ is twice differentiable, $\nabla^2h(\bar x)\in\R^{n\times n}$
denotes the Hessian of $h$ at $\bar x$.
Partial derivatives with respect to particular variables are represented in analogous fashion.

We use $\Gamma\colon\R^n\tto\R^m$ is order to express that $\Gamma$ is a so-called set-valued
mapping which assigns to each $x\in\R^n$ a (potentially empty) subset of $\R^m$.
For any such set-valued mapping, $\dom\Gamma:=\{x\in\R^n\,|\,\Gamma(x)\neq\emptyset\}$
and $\gph\Gamma:=\{(x,y)\in\R^n\times\R^m\,|\,y\in\Gamma(x)\}$ are the domain and the graph
of $\Gamma$, respectively.
We say that $\Gamma$ is inner semicontinuous at $(\bar x,\bar y)\in\gph\Gamma$ whenever for
each sequence $\{x^k\}_{k\in\N}\subset\R^n$ such that $x^k\to\bar x$, there exists a sequence
$\{y^k\}_{k\in\N}\subset\R^m$ such that $y^k\to\bar y$ which satisfies $y^k\in\Gamma(x^k)$
for all large enough $k\in\N$. We emphasize that whenever $\Gamma$ is inner semicontinuous at $(\bar x,\bar y)$,
then $\bar x$ is an interior point of $\dom\Gamma$.
At a fixed point $(\bar x,\bar y)\in\gph\Gamma$, $\Gamma$ is said to be calm whenever
there are constants $\varepsilon,\delta,L>0$ such that
\[
	\forall x\in\mathbb B_\varepsilon(\bar x),\,
	\forall y\in \Gamma(x)\cap\mathbb B_\delta(\bar y),\,
	\exists \tilde y\in\Gamma(\bar x)\colon\quad
	\norm{y-\tilde y}\leq L\norm{x-\bar x}.
\]
We note that whenever $H\colon\R^n\to\R^m$ is a single-valued function which is calm at $\bar x$,
then the set-valued mapping $x\mapsto\{H(x)\}$ is calm at $(\bar x,H(\bar x))$.
The converse is not necessarily true.

\subsection{Newton-differentiability}\label{sec:Newton_differentiability}

In order to construct numerical methods for the computational solution of nonsmooth equations,
one has to choose a suitable concept of generalized differentiation.
Typically, the idea of semismoothness is used for that purpose, see \cite{Mifflin1977,QiSun1993}.
Here, however, we strike a different path and exploit the notion of
Newton-differentiability, 
see \cite{HintermuellerItoKunisch2002,Ulbrich2002},
whose origins can be found in infinite-dimensional optimization. 
The latter has the inherent advantage that, in contrast to semismoothness, 
it is defined for non-Lipschitz functions and enjoys a natural calculus which 
follows the lines known from standard differentiability.
These beneficial features have been used recently in order to solve stationarity
conditions of complementarity-constrained optimization problems which can be reformulated
as systems of discontinuous equations, see \cite{HarderMehlitzWachsmuth2021}.

Let us start with the formal definition of Newton-differentiability before
presenting some essential calculus rules. Most of this material is taken from
the recent contribution \cite[Section~2.3]{HarderMehlitzWachsmuth2021}.
\begin{definition}\label{def:Newton_differentiability}
	Let $\calF\colon\R^p\to\R^q$ and $D_N\calF\colon\R^p\to\R^{q\times p}$
	be given mappings and let $M\subset\R^p$ be nonempty.
	We say that
	\begin{enumerate}
		\item $\calF$ is Newton-differentiable on $M$ with Newton-derivative $D_N\calF$ whenever
			for each $z\in M$, we have
			\[
				\calF(z+d)-\calF(z)-D_N\calF(z+d)d=\oo(\norm{d}),
			\]
		\item $\calF$ is Newton-differentiable on $M$ of order $\alpha\in(0,1]$ with Newton-derivative
			$D_N\calF$ whenever for each $z\in M$, we have
			\[
				\calF(z+d)-\calF(z)-D_N\calF(z+d)d=\OO(\norm{d}^{1+\alpha}).
			\]
	\end{enumerate}
\end{definition}

We note that the Newton-derivative of a mapping $\calF\colon\R^p\to\R^q$
is not necessarily uniquely determined if $\calF$ is
Newton-differentiable on some set $M\subset\R^p$.
It can be easily checked that any continuously differentiable function $\calF\colon\R^p\to\R^q$
is Newton-differentiable on $\R^p$ when $D_N\cal F:=\calF'$ is chosen. 
If, additionally, $\calF'$ is locally Lipschitz continuous, then the order of
Newton-differentiability is $1$.

In the recent paper \cite{HarderMehlitzWachsmuth2021}, the authors present several important
calculus rules for Newton-derivatives including a chain rule which we recall below, see
\cite[Lemma~2.11]{HarderMehlitzWachsmuth2021}.

\begin{lemma}\label{lem:chain_rule_Newton_derivative}
	Suppose that $\calF\colon\R^p\to\R^q$ is Newton-differentiable on $M\subset\R^p$
	with Newton-derivative $D_N\calF$, 
	and that $\calG\colon\R^q\to\R^s$ is Newton-differentiable on $\calF(M)$ 
	with Newton-derivative $D_N\calG$.
	Further, assume that $D_N\calF$ is bounded on a neighborhood of $M$,
	and that $D_N\calG$ is bounded on a neighborhood of $\calF(M)$.
	Then $\calG\circ\calF$ is Newton-differentiable on $M$ with
	Newton-derivative given by $z\mapsto D_N\calG(\calF(z))D_N\calF(z)$.
	If both $\calF$ and $\calG$ are Newton-differentiable of order $\alpha\in(0,1]$,
	then $\calG\circ\calF$ is Newton-differentiable of order $\alpha$ with the
	Newton-derivative given above.
\end{lemma}

Let us note that \cref{lem:chain_rule_Newton_derivative} directly gives a sum rule
which applies to the concept of Newton-differentiability. However, a direct proof 
via \cref{def:Newton_differentiability}
shows that a sum rule for Newton-differentiability holds without additional
boundedness assumptions on the Newton-derivatives of the involved functions.

Let us inspect some examples which will be of major interest in this paper.
\begin{example}\label{ex:Newton_differentiability}
	\leavevmode
	\begin{enumerate}
		\item\label{item:newton_differentiability_max} 
			The function $\max\colon\R^2\to\R$ is Newton-differentiable on $\R^2$ of order $1$
			with Newton-derivative given by
			\[
				\forall (a,b)\in\R^2\colon\quad
				D_N\max(\cdot,\cdot)(a,b)
				:=
				\begin{cases}
					(1,0)	&	a\geq b,\\
					(0,1)	&	a<b,
				\end{cases}
			\]
			see \cite[Example~2.8]{HarderMehlitzWachsmuth2021}.
		\item\label{item:newton_differentiability_norm}
			For arbitrary $p\in\N$, we investigate Newton-differentiability of the Euclidean norm
			$\norm{\cdot}\colon\R^p\to\R$. Noting that $\norm{\cdot}$ is continuously
			differentiable on $\R^p\setminus\{0\}$ with locally Lipschitzian derivative, 
			it is Newton-differentiable of order $1$ there. Let us consider the mapping
			$D_N\norm{\cdot}\colon\R^p\to\R^{1\times p}$ given by
			\begin{equation}\label{eq:Newton_differentiability_norm}
				\forall z\in\R^p\colon\quad
				D_N\norm{\cdot}(z)
				:=
				\begin{cases}
					\frac{z^\top}{\norm{z}}	&	z\neq 0,\\
					\frac{\sqrt p}{p}\mathtt e^\top			&	z=0.
				\end{cases}
			\end{equation}
			Here, $\mathtt e\in\R^p$ denotes the all-ones vector.
			At the origin, we have $\norm{d}-D_N\norm{\cdot}(d)d=0$ for each $d\in\R^p$,
			so we already obtain Newton-differentiability of $\norm{\cdot}$ on $\R^p$
			of order $1$.
			We also note that the precise value of $D_N\norm{\cdot}(0)$ is completely
			irrelevant for the validity of this property.
			However, the particular choice in \eqref{eq:Newton_differentiability_norm}
			will be beneficial later on.
		\item\label{item:newton_differentiability_FB}
			 Let us investigate Newton-differentiability of the aforementioned
			FB function
			$\varphi_\textup{FB}\colon\R^2\to\R$ given in \eqref{eq:FB_function}.
			Relying on the sum rule (note that $\varphi_\textup{FB}$ is the sum of the
			identity and $\norm{\cdot}$ in $\R^2$) and respecting our arguments from 
			\cref{ex:Newton_differentiability}\,\ref{item:newton_differentiability_norm}, we obtain that
			$\varphi_\textup{FB}$ is Newton-differentiable on $\R^2$ of order $1$ with
			Newton-derivative given by
			\[
				\forall (a,b)\in\R^2\colon\quad
				D_N\varphi_\textup{FB}(a,b)
				:=
				\begin{cases}
				\left(1+\frac{a}{\sqrt{a^2+b^2}},1+\frac{b}{\sqrt{a^2+b^2}}\right)	
					&	(a,b)\neq(0,0),\\
				\left(1+\frac{\sqrt 2}{2},1+\frac{\sqrt 2}{2}\right)	
					&	(a,b)=(0,0),
				\end{cases}
			\]
			where we made use of \eqref{eq:Newton_differentiability_norm}.
	\end{enumerate}
\end{example}

\section{A local nonsmooth Levenberg--Marquardt method}\label{sec:LM_local}

The paper \cite{DeLucaFacchineiKanzow1996} introduces a semismooth Newton-type method for square
complementarity systems whose globalization is based on the FB function.
An extension to LM-type methods (which even can handle inexact solutions
of the subproblems) can be found in \cite{FacchineiKanzow1997,Jiang1999}.
A theoretical and numerical comparison of these methods is provided in \cite{DeLucaFacchineiKanzow2000}.
An application of the semismooth LM method in the context of mixed complementarity systems
can be found in \cite{KanzowPetra2004}, and these ideas were applied in the context of semiinfinite optimization 
e.g.\ in \cite{MaWang2009,QiWuZhou2003}.

In this subsection, we aim to analyze the local convergence properties
of nonsmooth LM methods in a much broader context which covers
not only the applications within the aforementioned papers, but also some new ones
in bilevel optimization, see e.g.\ \cref{sec:optimistic_bpp} and our comments
in \cref{sec:conclusions}. Our approach via Newton-derivatives is slightly
different from the one in the literature which exploits semismoothness of the underlying
map of interest and, thus, particularly, local Lipschitz continuity.

Throughout the subsection, we assume that $\calF\colon\R^p\to\R^q$ is a given mapping with
$q>p$ and inherent nonsmooth structure. We aim to solve the (over-determined) nonlinear system of equations
\begin{equation}\label{eq:root_problem}
	\calF(z)=0.
\end{equation}
Classically, this can be done by minimizing the squared residual of a
first-order linearization associated with $\calF$.
The basic idea behind LM methods is now to add a classical square-norm
regularization term to this  optimization problem.
Let us consider a current iterate $z^k$ where $\calF$ is Newton-differentiable, and consider
the minimization of 
$d\mapsto\tfrac12\norm{\calF(z^k)+D_N\calF(z^k)d}^2+\tfrac{\nu_k}{2}\norm{d}^2$
which is a strictly convex function.
Above, we exploited the Newton-derivative in order to find a linearization of $\calF$
at $z^k$, and $\nu_k>0$ is a regularization parameter.
By means of the chain rule, a necessary optimality condition for this surrogate problem
is given by
\begin{equation}\label{eq:LM_equation}
	\left(D_N\calF(z^k)^\top D_N\calF(z^k)+\nu_k \mathbb I^p\right)d
	=
	-D_N\calF(z^k)^\top \calF(z^k).
\end{equation}
Observing that the matrix $D_N\calF(z^k)^\top D_N\calF(z^k)$ is at least positive semidefinite,
\eqref{eq:LM_equation} indeed possesses a uniquely determined solution $d^k$.

This motivates the formulation of \cref{alg:LocalLM} 
for the numerical treatment of \eqref{eq:root_problem}.
We assume that $D_N\calF\colon\R^p\to\R^{q\times p}$ is a given function which serves as a
Newton-derivative on a suitable subset of $\R^p$ which will be specified later.

\begin{algorithm}[Local nonsmooth Levenberg--Marquardt method]\leavevmode
	\label{alg:LocalLM}
	\begin{algorithmic}[1]
		\REQUIRE starting point $z^0\in\R^p$
		\STATE set $k := 0$
		\WHILE{a suitable termination criterion is violated at iteration $ k $}
		\STATE choose $\nu_k>0$ and compute $d^k$ as the uniquely determined
			solution of \eqref{eq:LM_equation}
		\STATE set $z^{k+1}:=z^k+d^k$ and $k:=k+1$
		\ENDWHILE
		\RETURN $z^k$
	\end{algorithmic}
\end{algorithm}

In the subsequent theorem, we present a local convergence result regarding
\cref{alg:LocalLM}.

\begin{theorem}\label{thm:local_convergence_LM}
	Let $M\subset\R^p$ be a set such that a solution $\bar z\in\R^p$ of
	\eqref{eq:root_problem} satisfies $\bar z\in M$.
	Assume that $\calF$ is Newton-differentiable on $M$ with
	Newton-derivative $D_N\calF\colon\R^p\to\R^{q\times p}$.
	We assume that there are $\delta>0$ and $C>0$
	such that, for all $z\in\mathbb B_\delta(\bar z)$, $D_N\calF(z)$ possesses
	full column rank $p$ such that
	\begin{equation}\label{eq:uniform_boundedness_eigenvalues}
		\lambda_\textup{min}(D_N\calF(z)^\top D_N\calF(z))\geq \frac{1}{C},\qquad
		\norm{D_N\calF(z)}\leq C.
	\end{equation}
	Then there exists $\varepsilon>0$ such that, for each starting point
	$z^0\in\mathbb B_\varepsilon(\bar z)$ and
	each null sequence $\{\nu_k\}_{k\in\N}\subset(0,\varepsilon)$,
	\cref{alg:LocalLM} terminates after finitely many steps
	or produces a sequence $\{z^k\}_{k\in\N}$ which converges
	to $\bar z$ superlinearly.
	Furthermore, if $\calF$ is even Newton-differentiable on $M$ of order $1$, and if
	$\nu_k\in\OO(\nnorm{\calF(z^k)})$ while $\calF$ is calm at $\bar z$,
	then the convergence is quadratic.	
\end{theorem}
\begin{proof}
	Due to the assumptions of the theorem, and by Newton-differentiability
	of $\calF$,
	we can choose $\varepsilon\in(0,\min(\delta,1/(4C))$ so small such that
	the following estimates hold true for all $d\in\mathbb B_\varepsilon(0)$
	and $\nu>0$:
	\begin{subequations}\label{eq:proof_local_LM}
		\begin{align}
			\label{eq:proof_local_LM_uniform_invertibility}
				\nnorm{(D_N\calF(\bar z+d)^\top D_N\calF(\bar z+d)+\nu \mathbb I^p)^{-1}}
				&\leq
				C,\\
			\label{eq:proof_local_LM_bounded_derivative}
				\norm{D_N\calF(\bar z+d)}
				&\leq
				C,\\
			\label{eq:proof_local_LM_Newton_differentiability}
				\norm{\calF(\bar z+d)-\calF(\bar z)-D_N\calF(\bar z+d)d}
				&\leq
				\frac{1}{4C^2}\norm{d}.
		\end{align}
	\end{subequations}
	Using $\calF(\bar z)=0$, for each $z^k\in\mathbb B_\varepsilon(\bar z)$ and
	$\nu_k\in(0,\varepsilon)$, we find
	\begin{equation}\label{eq:estimate_local_LM}
		\begin{aligned}
		\nnorm{z^{k+1}-\bar z}
		&=
		\nnorm{z^k-(D_N\calF(z^k)^\top D_N\calF(z^k)+\nu_k\mathbb I^p)^{-1}
		D_N\calF(z^k)^\top \calF(z^k)-\bar z}
		\\
		&\leq
		C\nnorm{D_N\calF(z^k)^\top(\calF(z^k)-\calF(\bar z)-D_N\calF(z^k)(z^k-\bar z))}
		+
		C\nu_k\nnorm{z^k-\bar z}
		\\
		&
		\leq
		C^2\nnorm{\calF(z^k)-\calF(\bar z)-D_N\calF(z^k)(z^k-\bar z)}+C\nu_k\nnorm{z^k-\bar z}
		\\
		&\leq
		\tfrac{1}{4}\nnorm{z^k-\bar z}+\tfrac{1}{4}\nnorm{z^k-\bar z}
		=
		\tfrac{1}{2}\nnorm{z^k-\bar z}.
		\end{aligned}
	\end{equation}
	Thus, we have $\nnorm{z^{k+1}-\bar z}\leq\frac12\nnorm{z^k-\bar z}$, i.e.,
	$z^{k+1}\in\mathbb B_{\varepsilon/2}(\bar z)$.
	Thus, if $z^0\in \mathbb B_\varepsilon(\bar z)$ and
	$\{\nu_k\}_{k\in\N}\subset(0,\varepsilon)$, we have $z^k\to\bar z$ in the case where
	\cref{alg:LocalLM} generates an infinite sequence.
	Furthermore, the definition of Newton-differentiability,
	\eqref{eq:estimate_local_LM}, and $\nu_k\to 0$ give
	\[
		\nnorm{z^{k+1}-\bar z}=\oo(\nnorm{z^k-\bar z}),
	\]
	i.e., the convergence $z^k\to\bar z$ is superlinear.
	
	Finally, assume that $\calF$ is Newton-differentiable of order $1$ and calm at $\bar z$, 
	and that $\nu_k\in\OO(\nnorm{\calF(z^k)})$.
	Then there is a constant $K>0$ such that the estimate
	\eqref{eq:estimate_local_LM} can be refined as
	\begin{align*}
		\nnorm{z^{k+1}-\bar z}
		&\leq
		C^2\nnorm{\calF(z^k)-\calF(\bar z)-D_N\calF(z^k)(z^k-\bar z)}+C\nu_k\nnorm{z^k-\bar z}
		\\
		&\leq
		\OO(\nnorm{z^k-\bar z}^2)+CK\nnorm{\calF(z^k)}\nnorm{z^k-\bar z}
		\\
		&=
		\OO(\nnorm{z^k-\bar z}^2),
	\end{align*}
	where we used $\calF(\bar z)=0$ and the calmness of $\calF$ at $\bar z$ in the last step.
\end{proof}

Let us briefly compare the assumptions of \cref{thm:local_convergence_LM}, which are used to
guarantee the superlinear or quadratic convergence of a sequence generated by \cref{alg:LocalLM},
with the ones exploited in the literature where nonsmooth LM methods are
considered from the viewpoint of semismoothness, see e.g.\ \cite[Section~2]{FacchineiKanzow1997}.
Therefore, we fix a solution $\bar z\in\R^p$ of \eqref{eq:root_problem}.
The full column rank assumption on the Newton-derivative, locally around $\bar z$, corresponds to so-called
BD-regularity of the point $\bar z$ which demands that all matrices within Bouligand's
generalized Jacobian at $\bar z$ (if available), see e.g.\ \cite[Section~2]{Qi1993}, possess full column rank.
We note that, by upper semicontinuity of Bouligand's generalized Jacobian, this full rank assumption
extends to a neighborhood of $\bar z$. Hence, these full rank assumptions are, essentially, 
of the same type although the one from \cref{thm:local_convergence_LM} is more general since it addresses
situations where the underlying map is allowed to be non-Lipschitz.
Second, \cref{thm:local_convergence_LM} assumes boundedness of the Newton-derivative, locally around $\bar z$.
In the context of semismooth LM methods, local boundedness of the generalized 
derivative holds inherently by construction of Bouligand's generalized Jacobian and local Lipschitzness of the
underlying mapping. 
Third, for superlinear convergence, the regularization parameters need to satisfy 
$\nu_k\in\OO(\nnorm{\calF(z^k)})$ in \cref{thm:local_convergence_LM}, and this assumption is also used in
\cite{FacchineiKanzow1997}. Similarly, since Newton-differentiability of order $1$ is a natural
counterpart of so-called strict semismoothness, the assumptions for quadratic convergence are
also the same.
Summarizing these impressions, \cref{alg:LocalLM} and \cref{thm:local_convergence_LM} generalize the
already existing theory on nonsmooth LM methods to a broader setting under
reasonable assumptions.
We also note that our analysis made no use of deeper underlying properties of the generalized derivative,
we mainly used its definition for our purposes.
However, it should be observed that a sophisticated choice of the Newton-derivative, 
which is not uniquely determined for a given map as mentioned in \cref{sec:Newton_differentiability},  
may lead to weaker assumptions in \cref{thm:local_convergence_LM} than a bad choice of it.
	Indeed, it even may happen that the assumptions of \cref{thm:local_convergence_LM} are valid for
	a particular choice of the Newton-derivative while they are violated for another one.
	Thus, when \cref{alg:LocalLM} (or a suitable globalization) is implemented, one has to keep
	in mind to choose the Newton-derivative in such a way that the assumptions of 
	\cref{thm:local_convergence_LM} are valid (if this is actually possible), 
	as the requirements of \cref{thm:local_convergence_LM} need to be satisfied for one
	particular choice of the Newton-derivative (and not for all possible choices).
	Clearly, this choice depends on structural properties of the nonsmooth mapping $\mathcal F$ under
	consideration and is, potentially, a laborious task as we will see in \cref{sec:LM},
	see \cite{HarderMehlitzWachsmuth2021} as well.	

The following corollary of \cref{thm:local_convergence_LM} shows that 
$\nnorm{\calF(z^{k+1})}\in\oo(\nnorm{\calF(z^k)})$ can be expected under reasonable assumptions, and
this will be of essential importance later on.
\begin{corollary}\label{cor:quotients_of_residuals}
	Under the assumptions of \cref{thm:local_convergence_LM} which guarantee the
	superlinear convergence of $\{z^k\}_{k\in\N}$ generated by
	\cref{alg:LocalLM}, we additionally have
	\[
		\frac{\nnorm{\calF(z^{k+1})}}{\norm{\calF(z^k)}}\to 0
	\]
	provided $\calF$ is calm at $\bar z$.
\end{corollary}
\begin{proof}
	We choose $\varepsilon>0$ as in the proof of \cref{thm:local_convergence_LM} and observe
	$\{z^k\}_{k\in\N}\subset\mathbb B_\varepsilon(\bar z)$.
	Exploiting the Newton-differentiability of $\calF$ and $\calF(\bar z)=0$, we have
	\[
		\calF(z^k)=D_N\calF(z^k)(z^k-\bar z)+\oo(\nnorm{z^k-\bar z}),
	\]
	and some transformations give, for sufficiently large $k\in\N$ and by
	boundedness of the sequence $\{D_N\calF(z^k)\}_{k\in\N}$,
	\begin{align*}
		\nnorm{z^k-\bar z}
		\leq
		\nnorm{(D_N\calF(z^k)^\top D_N\calF(z^k))^{-1}D_N\calF(z^k)^\top \calF(z^k)}
		+
		\tfrac12\nnorm{z^k-\bar z},
	\end{align*}
	i.e.,
	\begin{align*}
		\tfrac12\nnorm{z^k-\bar z}
		\leq
		\nnorm{(D_N\calF(z^k)^\top D_N\calF(z^k))^{-1}}\nnorm{D_N\calF(z^k)}\nnorm{\calF(z^k)}
		\leq
		C^2\nnorm{\calF(z^k)}
	\end{align*}
	due to \eqref{eq:uniform_boundedness_eigenvalues}.
	Thus, we find
	\begin{align*}
		\frac{\nnorm{\calF(z^{k+1})}}{\nnorm{\calF(z^k)}}
		&\leq
		\frac{2C^2\nnorm{\calF(z^{k+1})-\calF(\bar z)}}{\nnorm{z^k-\bar z}}
		\leq
		\frac{2C^2L\nnorm{z^{k+1}-\bar z}}{\nnorm{z^k-\bar z}}
		\leq
		\frac{2C^2L\,\oo(\nnorm{z^k-\bar z})}{\nnorm{z^k-\bar z}}
		\to
		0
	\end{align*}
	as $k\to\infty$, where $L>0$ is a local calmness constant of $\calF$ at $\bar z$.
\end{proof}

	In \cite[Section~3.2]{HarderMehlitzWachsmuth2021}, a function is
	constructed which possesses the following properties:
	\begin{enumerate}
		\item[(i)] it is Newton-differentiable on its set of roots with globally
			bounded and nonvanishing Newton-derivative,
		\item[(ii)] it is discontinuous in each open neighborhood of its set of roots, and
		\item[(iii)] it is calm at each point from its set of roots.
	\end{enumerate}
	This function is then used to construct a nonsmooth Newton-type method for the
	computational solution of stationarity systems associated with
	complementarity-constrained optimization problems.
	We note that it crucially violates the standard requirement of local Lipschitzness
	which is typically exploited in the context of nonsmooth Newton-type methods.
	Similarly as in \cite{HarderMehlitzWachsmuth2021},
	the analysis in \cref{thm:local_convergence_LM} and \cref{cor:quotients_of_residuals} 
	only requires calmness of the mapping $\calF$ (as well as some other natural assumptions) 
	in order to get local fast convergence of a nonsmooth LM method. 
Thus, the ideas from \cite{HarderMehlitzWachsmuth2021} can be
carried over to the situation where over-determined stationarity systems
of complementarity-constrained optimization problems need to be solved
(such systems would, exemplary, arise when applying the theory from \cref{sec:obpp_modeling}
to bilevel optimization problems with additional complementarity constraints at the
upper-level stage or to the so-called combined reformulation of bilevel optimization
problems which makes use of the so-called value function and Karush--Kuhn--Tucker reformulation
at the same time, see \cite{YeZhu2010}).

For the globalization of \cref{alg:LocalLM}, 
one typically needs to impose additional assumptions like the
smoothness of $z\mapsto\norm{\calF(z)}^2$.
In the next section, we address the prototypical setting 
of mixed nonlinear complementarity systems and carve out
two classical globalization strategies.

\section{A global nonsmooth Levenberg--Marquardt method for mixed nonlinear complementarity systems}\label{sec:LM}

For continuously differentiable functions
$\calH\colon\R^{p_1}\times\R^{p_2}\to\R^{q_1}$ and $\calG\colon\R^{p_1}\times\R^{p_2}\to\R^{p_2}$,
we aim to solve the mixed nonlinear complementarity system
\begin{equation}\label{eq:NLCS}\tag{MNLCS}
	\calH(w,\xi)=0,\qquad \calG(w,\xi)\leq 0,\,\xi\geq 0,\,\calG(w,\xi)^\top\xi=0
\end{equation}
where, in the application we have in mind, $q_1>p_1$ holds true,
see \cref{sec:optimistic_bpp}.
A comprehensive overview of available theoretical and numerical aspects 
as well as applications of mixed nonlinear complementarity problems can 
be found in the monograph \cite{FacchineiPang2003}.
Typical solution approaches are based on nonsmooth Newton-type methods,
see \cite{Kanzow2000}, smoothing methods, see \cite{ChenMangasarian1996,LiFukushima2000},
active-set and interior-point methods, see \cite{DarynkaIzmailovSolodov2005,MonteiroPang1996},
as well as penalty methods, see \cite{HuangWang2012}.
Let us note that the classical formulation of mixed nonlinear complementarity systems as considered in
\cite{ChenMangasarian1996,DarynkaIzmailovSolodov2005,Kanzow2000,LiFukushima2000,MonteiroPang1996}
is essentially equivalent to the one we are suggesting in \eqref{eq:NLCS},
see \cite[Section~4]{ChenMangasarian1996} and \cite[Section~1]{DarynkaIzmailovSolodov2005} 
for arguments which reveal this equivalence.
However, in \cref{sec:optimistic_bpp}, it will be beneficial to work with the model \eqref{eq:NLCS}
to preserve structural properties.

For later use, we introduce some index sets depending on a pair
$(w,\xi)\in\R^{p_1}\times\R^{p_2}$:
\begin{align*}
	I^0(w,\xi)&:=\{i\in\{1,\ldots,p_2\}\,|\,\calG_i(w,\xi)=0\},\\
	I^-(w,\xi)&:=\{i\in\{1,\ldots,p_2\}\,|\,\calG_i(w,\xi)<0\},\\
	I^+(w,\xi)&:=\{i\in I^0(w,\xi)\,|\,\xi_i>0\},\\
	I^{00}(w,\xi)&:=\{i\in I^0(w,\xi)\,|\,\xi_i=0\}.
\end{align*}
Above, $\mathcal G_1,\dots,\mathcal G_{p_2}\colon\R^{p_1}\times\R^{p_2}\to\R$ 
are the component functions of $\mathcal G$.

Based on NCP-functions, see \cref{sec:introduction}, the complementarity condition
\begin{equation}\label{eq:pure_complementarity_constraints}
	\calG(w,\xi)\leq 0,\,\xi\geq 0,\,\calG(w,\xi)^\top\xi=0
\end{equation}
can be restated in form of a (nonsmooth) equality constraint.
In this paper, we will focus on two popular choices, namely the maximum function and
the famous FB function.
Clearly, \eqref{eq:pure_complementarity_constraints} is equivalent to
\[
	\max(\calG(w,\xi),-\xi)=0,
\]
where $\max\colon\R^{p_2}\times\R^{p_2}\to\R^{p_2}$ is exploited to express the componentwise maximum.
Using $\phi_\textup{FB}\colon\R^{p_2}\times\R^{p_2}\to\R^{p_2}$, given by
\[
	\forall \xi^1,\xi^2\in\R^{p_2}\colon\quad
	\phi_\textup{FB}(\xi^1,\xi^2)
	:=
	\begin{pmatrix}
		\varphi_\textup{FB}(\xi^1_1,\xi^2_1)
		\\
		\vdots
		\\
		\varphi_\textup{FB}(\xi^1_{p_2},\xi^2_{p_2})
	\end{pmatrix}
\]
where $\varphi_\textup{FB}\colon\R^2\to\R$ is the FB function 
defined in \eqref{eq:FB_function},
\eqref{eq:pure_complementarity_constraints} is also equivalent to
\[
	\phi_\textup{FB}(\calG(w,\xi),-\xi)=0.
\]
This motivates the consideration of the special residual mappings
$\calF_\textup{max},\calF_\textup{FB}\colon\R^{p_1}\times\R^{p_2}\to\R^{q_1}\times\R^{p_2}$
given by
\begin{equation}\label{eq:special_residual_complementarity}
	\begin{aligned}
	&\forall w\in\R^{p_1},\,\forall \xi\in\R^{p_2}\colon\quad&
	\calF_{\textup{max}}(w,\xi)
	&:=
	\begin{bmatrix}
		\calH(w,\xi)\\
		\max(\calG(w,\xi),-\xi)
	\end{bmatrix},
	&
	\\
	&&
	\calF_\textup{FB}(w,\xi)
	&:=
	\begin{bmatrix}
		\calH(w,\xi)\\
		\phi_\textup{FB}(\calG(w,\xi),-\xi)
	\end{bmatrix},
	&
	\end{aligned}
\end{equation}
since precisely their roots are the solutions of \eqref{eq:NLCS}.
Let us note that, in a certain sense, these residuals are equivalent.
This can be distilled from \cite[Lemma~3.1]{Tseng1996}.
\begin{lemma}\label{lem:equivalence_of_residuals}
	There exist constants $K_1,K_2>0$ such that
	\[
		\forall w\in\R^{p_1},\,\forall \xi\in\R^{p_2}\colon\quad
		K_1\norm{\calF_{\textup{max}}(w,\xi)}
		\leq
		\nnorm{\calF_\textup{FB}(w,\xi)}
		\leq 
		K_2\norm{\calF_{\textup{max}}(w,\xi)}.
	\]
\end{lemma}

Due to the inherent nonsmoothness of $\calF_{\textup{max}}$ and $\calF_\textup{FB}$
from \eqref{eq:special_residual_complementarity}, we may now
apply \cref{alg:LocalLM} in order to solve $\calF_{\textup{max}}(w,\xi)=0$ 
or $\calF_\textup{FB}(w,\xi)=0$.
Note that we have $p:=p_1+p_2$ and $q:=q_1+p_2$ in the context of \cref{sec:LM_local}.
In the literature, the approach via $\calF_\textup{FB}$ is quite popular due to the
well-known observation that the mapping $\Psi_\textup{FB}\colon\R^{p_1}\times\R^{p_2}\to\R$
given by
\[
	\forall w\in\R^{p_1},\,\forall \xi\in\R^{p_2}\colon\quad
	\Psi_\textup{FB}(w,\xi)
	:=
	\frac12\nnorm{\calF_\textup{FB}(w,\xi)}^2
\]
is continuously differentiable, allowing for a globalization of suitable local
solution methods. This is due to the fact that the squared FB function
is continuously differentiable, see e.g.\ \cite[Lemma~3.4]{FacchineiSoares1997}. 
However, it has been reported in
\cite{DeLucaFacchineiKanzow2000} in the context of square nonlinear complementarity
systems that a mixed approach combining both residuals
from \eqref{eq:special_residual_complementarity} might be beneficial since the assumptions
for local fast convergence could be weaker while the local convergence is faster.

\subsection{The mixed approach via both residual functions}\label{sec:global_LM_max}

First, let us focus on this mixed approach
where the central LM search direction
will be computed via $\calF_{\textup{max}}$ given in \eqref{eq:special_residual_complementarity}
in order to end up with less complicated Newton-derivatives
in the LM system. Still, we exploit the smooth function $\Psi_\textup{FB}$
for globalization purposes.
Later, in \cref{sec:global_LM_FB}, we briefly comment on the method which exclusively uses $\calF_\textup{FB}$.

Below, we provide formulas for a Newton-derivative of $\calF_{\textup{max}}$ and the gradient of
$\Psi_\textup{FB}$ for later use.
\begin{lemma}\label{lem:Newton_derivative_of_complementarity_system}
	On each bounded subset of $\R^{p_1}\times\R^{p_2}$, the mapping $\calF_{\textup{max}}$ is
	Newton-differentiable with Newton-derivative given by
	\[
		(w,\xi)
		\mapsto
		\begin{bmatrix}
			\calH'_w(w,\xi)&\calH'_\xi(w,\xi)\\
			\mathbb D(I^\geq(w,\xi))\calG'_w(w,\xi)&\mathbb D(I^{\geq}(w,\xi))\calG'_\xi(w,\xi)-\mathbb D(I^<(w,\xi))
		\end{bmatrix},
	\]
	and this mapping, again, is bounded on bounded sets.
	Whenever the derivatives of $\calG$ and $\calH$ are locally Lipschitzian,
	the order of Newton-differentiability is $1$.
	Above, the index sets $I^\geq(w,\xi)$ and $I^<(w,\xi)$ are given by
	\begin{align*}
		I^\geq(w,\xi)
		&:=
		\{i\in\{1,\ldots,p_2\}\,|\,\calG_i(w,\xi)\geq-\xi_i\},\\
		I^<(w,\xi)
		&:=
		\{i\in\{1,\ldots,p_2\}\,|\,\calG_i(w,\xi)<-\xi_i\},
	\end{align*}
	and, for each $I\subset\{1,\ldots,p_2\}$, $\mathbb D(I)\in\R^{p_2\times p_2}$
	is the diagonal matrix given by
	\[
		\forall i,j\in\{1,\ldots,p_2\}\colon\quad
		\mathbb D(I)_{ij}
		:=
		\begin{cases}
			1	&i=j,\,i\in I,\\
			0	&\text{otherwise.}
		\end{cases}
	\]
\end{lemma}
\begin{proof}
	This follows easily from \cref{lem:chain_rule_Newton_derivative},
	\cref{ex:Newton_differentiability}\,\ref{item:newton_differentiability_max},
	and our comments right after \cref{def:Newton_differentiability}.
\end{proof}

Throughout the section, we will denote the Newton-derivative of $\calF_\textup{max}$ which
has been characterized in \cref{lem:Newton_derivative_of_complementarity_system} by
$D_N\calF_{\textup{max}}$.

The next result follows simply by computing the derivative of the squared FB
function and using the standard chain rule.
\begin{lemma}\label{lem:gradient_psi_FB}
	At each point $(w,\xi)\in\R^{p_1}\times\R^{p_2}$, $\Psi_\textup{FB}$
	is continuously differentiable, and its gradient is given by
	\[
		\nabla\Psi_\textup{FB}(w,\xi)
		=
		\begin{bmatrix}
			\calH'_w(w,\xi)&\calH'_\xi(w,\xi)\\
			\widetilde{\mathbb D}_\calG(w,\xi)\calG'_w(w,\xi)&\widetilde{\mathbb D}_\calG(w,\xi)\calG'_\xi(w,\xi)-\widetilde{\mathbb D}_\xi(w,\xi)
		\end{bmatrix}^\top
		\calF_\textup{FB}(w,\xi)
	\]
	where
	\[
		\widetilde{\mathbb D}_\calG(w,\xi):=\diag(v_a(w,\xi)),
		\qquad
		\widetilde{\mathbb D}_\xi(w,\xi):=\diag(v_b(w,\xi)).
	\]
	Above, $v_a(w,\xi),v_b(w,\xi)\in\R^{p_2}$ are the vectors given
	by
	\begin{equation}\label{eq:surrogate_vectors_FB}
	\begin{aligned}
		\forall i\in\{1,\ldots,p_2\}\colon\quad
		(v_a(w,\xi))_i
		&:=
		\begin{cases}
			1+\frac{\calG_i(w,\xi)}{\sqrt{\calG_i^2(w,\xi)+\xi_i^2}}	&	i\notin I^{00}(w,\xi),\\
			1+\frac{\sqrt 2}{2}									&	i\in I^{00}(w,\xi),
		\end{cases}
		\\
		(v_b(w,\xi))_i
		&:=
		\begin{cases}
			1-\frac{\xi_i}{\sqrt{\calG_i^2(w,\xi)+\xi_i^2}}	&	i\notin I^{00}(w,\xi),\\
			1+\frac{\sqrt 2}{2}							&	i\in I^{00}(w,\xi).
		\end{cases}
	\end{aligned}
	\end{equation}
\end{lemma}

For local superlinear convergence of our method of interest, we need the following result.
\begin{lemma}\label{lem:local_uniform_invertibility}
	Fix a solution $(\bar w,\bar\xi)\in\R^{p_1}\times\R^{p_2}$ of \eqref{eq:NLCS} and
	assume that for each index set $I\subset I^{00}(\bar w,\bar\xi)$, the matrix
	\[
		\begin{bmatrix}
			\calH'_w(\bar w,\bar\xi)	&	\calH'_\xi(\bar w,\bar \xi)\\
			\calG'_w(\bar w,\bar\xi)_{I^+(\bar w,\bar\xi)\cup I}&\calG'_\xi(\bar w,\bar\xi)_{I^+(\bar w,\bar\xi)\cup I}\\
			\mathbb O	&-\mathbb I^{p_2}_{I^-(\bar w,\bar\xi)\cup I_\textup{c}}
		\end{bmatrix},
	\]
	where we used $I_\textup{c}:=I^{00}(\bar w,\bar\xi)\setminus I$,
	possesses linear independent columns.
	Then the matrices $\mathcal M(w,\xi)^\top\mathcal M(w,\xi)$, where
	\[
		\mathcal M(w,\xi)
		:=
		\begin{bmatrix}
			\calH'_w(w,\xi)&\calH'_\xi(w,\xi)\\
			\mathbb D(I^\geq(w,\xi))\calG'_w(w,\xi)&\mathbb D(I^\geq(w,\xi))\calG'_\xi(w,\xi)-\mathbb D(I^<(w,\xi))
		\end{bmatrix},
	\]
	are uniformly positive definite in a neighborhood of $(\bar w,\bar\xi)$.
\end{lemma}
\begin{proof}
	Suppose that the assertion is false.
	Then there exist sequences $\{w^k\}_{k\in\N}\subset\R^{p_1}$, $\{\xi^k\}_{k\in\N}\subset\R^{p_2}$,
	$\{(d^k_w,d^k_\xi)\}_{k\in\N}\subset\R^{p_1}\times\R^{p_2}$, and $\{\eta_k\}_{k\in\N}\subset[0,\infty)$
	such that $w^k\to\bar w$, $\xi^k\to\bar\xi$, $\eta_k\to 0$, and, for each $k\in\N$,
	$\nnorm{d^k_w}+\nnorm{d^k_\xi}=1$ as well as
	\[
		\eta_k
		=
		\begin{bmatrix}d^k_w\\d^k_\xi\end{bmatrix}^\top
		\mathcal M(w^k,\xi^k)^\top
		\mathcal M(w^k,\xi^k)
		\begin{bmatrix}d^k_w\\d^k_\xi\end{bmatrix}
		.
	\]
	For brevity of notation, we make use of the abbreviations
	\begin{equation}\label{eq:some_abbrevs}
		\begin{aligned}
			\calH'_w(k)&:=\calH'_w(w^k,\xi^k),&
			\qquad
			\calH'_\xi(k)&:=\calH'_\xi(w^k,\xi^k),&
			\\
			(\calG_i)'_w(k)&:=(\calG_i)'_w(w^k,\xi^k),&
			\qquad
			(\calG_i)'_\xi(k)&:=(\calG_i)'_\xi(w^k,\xi^k)&
		\end{aligned}
	\end{equation}
	for all $i\in\{1,\ldots,p_2\}$.
	Due to $\mathbb D(I^\geq(w^k,\xi^k))\mathbb D(I^<(w^k,\xi^k))=\mathbb O$, the above gives
	\begin{equation}\label{eq:surrogate_system_uniform_invertibility}
		\begin{aligned}
		\eta_k
		=
		\nnorm{\calH'_w(k)d^k_w}^2
		&+
		2(\calH'_w(k)d^k_w)^\top(\calH'_\xi(k)d^k_\xi)
		+
		\nnorm{\calH'_\xi(k)d^k_\xi}^2
		\\
		&+
		\sum_{i\in I^\geq(w^k,\xi^k)}\bigl((\calG_i)'_w(k)d^k_w+(\calG_i)'_\xi(k)d^k_\xi\bigr)^2
		+
		\sum_{i\in I^<(w^k,\xi^k)}(d^k_\xi)^2_i.
		\end{aligned}
	\end{equation}
	By continuity of $\calG$, each index $i\in I^+(\bar w,\bar \xi)$ lies in
	$I^\geq(w^k,\xi^k)$ for sufficiently large $k\in\N$.
	Furthermore, any $i\in I^-(\bar w,\bar\xi)$ lies in $I^<(w^k,\xi^k)$ for sufficiently
	large $k\in\N$.
	For $i\in I^{00}(\bar w,\bar\xi)$, two scenarios are possible.
	Either there is an infinite subset $K^i\subset\N$ such that
	$i\in I^\geq(w^k,\xi^k)$ for all $k\in K^i$,
	or
	$i\in I^<(w^k,\xi^k)$ holds for all large enough $k\in\N$.
	Anyhow, since there are only finitely many indices in $\{1,\ldots,p_2\}$, we may
	choose an infinite subset $K\subset\N$ as well as an index
	set $I\subset I^{00}(\bar w,\bar \xi)$ such that
	$I^\geq(w^k,\xi^k)=I^+(\bar w,\bar\xi)\cup I$ and
	$I^<(w^k,\xi^k)=I^-(\bar w,\bar\xi)\cup I_\textup{c}$ is valid for each $k\in K$,
	where $I_\textup c:=I^{00}(\bar w,\bar\xi)\setminus I$.
	Hence, for each $k\in K$, \eqref{eq:surrogate_system_uniform_invertibility}
	is equivalent to
	\begin{equation}\label{eq:surrogate_system_uniform_invertibility_II}
		\begin{aligned}
		\eta_k
		=
		\nnorm{\calH'_w(k)d^k_w}^2
		&+
		2(\calH'_w(k)d^k_w)^\top(\calH'_\xi(k)d^k_\xi)
		+
		\nnorm{\calH'_\xi(k)d^k_\xi}^2
		\\
		&+
		\sum_{i\in I^+(\bar w,\bar\xi)\cup I}\bigl((\calG_i)'_w(k)d^k_w+(\calG_i)'_\xi(k)d^k_\xi\bigr)^2
		+
		\sum_{i\in I^-(\bar w,\bar\xi)\cup I_{\textup{c}}}(d^k_\xi)^2_i.
		\end{aligned}
	\end{equation}
	Clearly, along a subsubsequence (without relabeling), $\{(d^k_w,d^k_\xi)\}_{k\in K}$
	converges to some $(d_w,d_\xi)\in\R^{p_1}\times\R^{p_2}$ such that
	$\nnorm{d_w}+\nnorm{d_\xi}=1$.
	Thus, taking the limit $k\to_K\infty$ in
	\eqref{eq:surrogate_system_uniform_invertibility_II} gives
	\begin{align*}
		0
		=
		\nnorm{\calH'_w(\bar w,\bar\xi)d_w}^2
		&+
		2(\calH'_w(\bar w,\bar\xi)d_w)^\top(\calH'_\xi(\bar w,\bar\xi)d_\xi)
		+
		\nnorm{\calH'_\xi(\bar w,\bar\xi)d_\xi}^2
		\\
		&+
		\sum_{i\in I^+(\bar w,\bar\xi)\cup I}\bigl((\calG_i)'_w(\bar w,\bar\xi)d_w+(\calG_i)'_\xi(\bar w,\bar\xi) d_\xi\bigr)^2
		+
		\sum_{i\in I^-(\bar w,\bar\xi)\cup I_\textup{c}}(d_\xi)^2_i.
	\end{align*}
	This implies that the matrix $\widetilde{\mathcal M}_I(\bar w,\bar\xi)^\top\widetilde{\mathcal M}_I(\bar w,\bar \xi)$, where
	\[
		\widetilde{\mathcal M}_I(\bar w,\bar \xi)
		:=
		\begin{bmatrix}
			\calH'_w(\bar w,\bar\xi)&\calH'_\xi(\bar w,\bar\xi)\\
			\mathbb D(I^+(\bar w,\bar\xi)\cup I)\calG'_w(\bar w,\bar\xi)
			&
			\mathbb D(I^+(\bar w,\bar\xi)\cup I)\calG'_\xi(\bar w,\bar\xi)-\mathbb D(I^-(\bar w,\bar\xi)\cup I_\textup{c})
		\end{bmatrix},
	\]
	which is naturally positive semidefinite by construction,
	is not positive definite. Thus, $\widetilde{\mathcal M}_I(\bar w,\bar\xi)$	
	cannot possess full column rank, contradicting the lemma's assumptions.
\end{proof}

The qualification condition postulated in \cref{lem:local_uniform_invertibility}
actually corresponds to the linear independence of the columns of all
elements of Bouligand's generalized Jacobian of $\calF_{\textup{max}}$ at $(\bar w,\bar\xi)$,
i.e., these assumptions recover the BD-regularity condition from the
literature, and the latter is well established in the context of solution algorithms
for nonsmooth systems. 
Let us also mention that it can be easily checked by means of simple examples that
full column rank of the Newton-derivative of $\calF_{\textup{max}}$ at $(\bar w,\bar\xi)$, as constructed
in \cref{lem:Newton_derivative_of_complementarity_system}, does, in general,
not guarantee the uniform positive definiteness which has been shown under
the assumptions of \cref{lem:local_uniform_invertibility}.

We note that each solution of \eqref{eq:NLCS} is a global minimizer of $\Psi_\textup{FB}$
and, thus, a stationary point of this function.
The converse statement is not likely to hold.
Even for quadratic systems, one needs strong additional assumptions in order to
obtain such a result. 

Next, we present the method of our interest in \cref{alg:GlobalLM}.
For brevity of notation, we introduce
\[
	z:=\begin{bmatrix}w\\ \xi\end{bmatrix},\qquad
	d:=\begin{bmatrix}\delta w\\ \delta\xi\end{bmatrix},
\]
and similarly, we define $z^k$ and $d^k$.
Recall that $D_N\calF_{\textup{max}}(z)$ denotes the
Newton-derivative of $\calF_{\textup{max}}$ at $z$ characterized
in \cref{lem:Newton_derivative_of_complementarity_system}.

\begin{algorithm}[Global nonsmooth Levenberg--Marquardt method for mixed nonlinear complementarity systems]\leavevmode
	\label{alg:GlobalLM}
	\setcounter{ALC@unique}{0}
	\begin{algorithmic}[1]
		\REQUIRE starting point $z^0\in\R^{p_1}\times\R^{p_2}$,
			parameters $\kappa\in(0,1)$, $\tau_\textup{abs}>0$, $\sigma,\beta\in(0,1)$,
			$\rho_1,\rho_2>0$, and $\gamma_1,\gamma_2>0$
		\STATE set $k := 0$
		\WHILE{$\nnorm{\calF_\textup{FB}(z^k)}\geq\tau_\textup{abs}$}
		\STATE set $\nu_k:=\min(\gamma_1,\gamma_2\nnorm{\calF_\textup{FB}(z^k)})$ and compute
			$d^k$ as the uniquely determined
			solution of
			\[
				\bigl(D_N\calF_{\textup{max}}(z^k)^\top D_N\calF_{\textup{max}}(z^k)+\nu_k\mathbb I^p\bigr)d=-D_N\calF_{\textup{max}}(z^k)^\top \calF_{\textup{max}}(z^k)
			\]
		\IF{$\Psi_\textup{FB}(z^k+d^k)\leq \kappa\Psi_\textup{FB}(z^k)$}\label{step:ratio_test}
			\STATE set $z^{k+1}:=z^k+d^k$
		\ELSE
			\IF{$\Psi'_\textup{FB}(z^k)d^k>-\rho_1\nnorm{\nabla\Psi_\textup{FB}(z^k)}\nnorm{d^k}$ or $\nnorm{d^k}<\rho_2$}
			\label{step:angle_test_+}
				\STATE set $d^k:=-\nabla\Psi_\textup{FB}(z^k)$
			\ENDIF
			\STATE set $\alpha_k:=\beta^{i_k}$ where $i_k\in\N$ is the smallest positive integer
				such that
				\[
					\Psi_\textup{FB}(z^k+\beta^{i_k}d^k)
					\leq
					\Psi_\textup{FB}(z^k)
					+
					\beta^{i_k}\sigma\,\Psi_\textup{FB}'(z^k)d^k
				\]
			\STATE set $z^{k+1}:=z^k+\alpha_kd^k$
		\ENDIF
		\STATE set $k:=k+1$
		\ENDWHILE
		\RETURN $z^k$
	\end{algorithmic}
\end{algorithm}

We now present the central convergence result associated with \cref{alg:GlobalLM}.
Its proof is similar to the one of \cite[Theorem~5.2]{HarderMehlitzWachsmuth2021} but
included for the purpose of completeness.
\begin{theorem}\label{thm:global_convergence_LM}
	Let $\{z^k\}_{k\in\N}$ be a sequence generated by \cref{alg:GlobalLM}.
	\begin{enumerate}
		\item If $\Psi_\textup{FB}(z^k+d^k)\leq \kappa\Psi_\textup{FB}(z^k)$ holds infinitely many
			times in \cref{step:ratio_test},
			then $\{\Psi_\textup{FB}(z^k)\}_{k\in\N}$ is a null sequence and
			each accumulation point of $\{z^k\}_{k\in\N}$ is a solution of
			\eqref{eq:NLCS}.
		\item Each accumulation point of $\{z^k\}_{k\in\N}$ is a stationary point
			of $\Psi_\textup{FB}$.
		\item If an accumulation point $\bar z$ of $\{z^k\}_{k\in\N}$ satisfies
			the assumptions of \cref{lem:local_uniform_invertibility}, then
			the whole sequence $\{z^k\}_{k\in\N}$ converges to $\bar z$
			superlinearly.
			If the derivatives of $\calG$ and $\calH$ are locally Lipschitz continuous
			functions, then the convergence is even quadratic.
	\end{enumerate}
\end{theorem}
\begin{proof}
	\begin{enumerate}
		\item We note that \cref{alg:GlobalLM} is a descent method with respect to
			$\Psi_\textup{FB}$ which is bounded from below by $0$.
			Thus, the assumptions guarantee $\Psi_\textup{FB}(z^k)\to 0$.
			Noting that $\Psi_\textup{FB}$ is continuous, each accumulation point
			$\bar z$ of $\{z^k\}_{k\in\N}$ satisfies $\Psi_\textup{FB}(\bar z)=0$
			in this situation,
			giving $\calF_\textup{FB}(\bar z)=0$, and this means that $\bar z$ solves
			\eqref{eq:NLCS}.
		\item If the assumptions of the first statement hold, the assertion is clear.
			Thus, we may assume that, along the tail of the sequence,
			$\Psi_\textup{FB}(z^k+d^k)> \kappa\Psi_\textup{FB}(z^k)$ is valid.
			Assume without loss of generality that $\{z^k\}_{k\in \N}$ converges
			to some point $\bar z$.
			We proceed by a distinction of cases.\\
			If $\liminf_{k\to\infty}\nnorm{d^k}=0$,
			\cref{alg:GlobalLM} automatically gives $d^k=-\nabla\Psi_\textup{FB}(z^k)$ along a
			subsequence (without relabeling). Taking the limit along this subsequence
			gives $\nabla\Psi_\textup{FB}(\bar z)=0$ by continuity of $\nabla\Psi_\textup{FB}$.\\
			Next, consider the case $\liminf_{k\to\infty}\alpha_k\nnorm{d^k}>0$.
			Noting that $\{\Psi_\textup{FB}(z^k)\}_{k\in\N}$ is monotonically decreasing and
			bounded from below, this sequence converges.
			Furthermore, for all large enough iterations $k\in\N$, 
			the choice of the step size and the fact that the search direction is a descent
			direction for $\Psi_\textup{FB}$ guarantee that
			\[
				0
				\leq
				-\alpha_k\sigma\Psi_\textup{FB}'(z^k)d^k
				\leq
				\Psi_\textup{FB}(z^k)-\Psi_\textup{FB}(z^{k+1}).
			\]
			Thus, since $\{\Psi_\textup{FB}(z^k)\}_{k\in\N}$ is a Cauchy sequence,
			we find $\alpha_k\Psi_\textup{FB}'(z^k)d^k\to 0$.
			Noting that, by construction of \cref{alg:GlobalLM}, each search direction
			passes the angle test 
			$\Psi_\textup{FB}'(z^k)d^k\leq-\rho_1\nnorm{\nabla\Psi_\textup{FB}(z^k)}\nnorm{d^k}$
			for large enough $k\in\N$, the estimate
			\[
				\alpha_k\Psi_\textup{FB}'(z^k)d^k
				\leq 0
				\leq
				\rho_1\alpha_k\nnorm{\nabla\Psi_\textup{FB}(z^k)}\nnorm{d^k}
				\leq
				-\alpha_k\Psi_\textup{FB}'(z^k)d^k
			\]
			follows. Taking the limit $k\to\infty$ and noting that 
			$\{\alpha_k\nnorm{d^k}\}_{k\in\N}$ is bounded away from zero,
			$\nabla\Psi_\textup{FB}(z^k)\to 0$ follows, so 
			$\nabla\Psi_\textup{FB}(\bar z)=0$ is obtained 
			from continuity of $\nabla\Psi_\textup{FB}$.\\
			Finally, we assume that $\liminf_{k\to\infty}\nnorm{d^k}>0$ and
			$\liminf_{k\to\infty}\alpha_k\nnorm{d^k}=0$.
			For simplicity of notation, let $\alpha_k\nnorm{d^k}\to 0$.
			Then we have $\alpha_k\to 0$ by assumption.
			Particularly, for large enough $k\in\N$, the step size candidate
			$\beta^{-1}\alpha_k$ is rejected, i.e.,
			\[
					\Psi_\textup{FB}(z^k+\beta^{-1}\alpha_kd^k)
					>
					\Psi_\textup{FB}(z^k)
					+
					\beta^{-1}\alpha_k\sigma\,\Psi_\textup{FB}'(z^k)d^k.
			\]
			The smoothness of $\Psi_\textup{FB}$ allows for the application of the mean value
			theorem in order to find $\theta_k\in(0,\beta^{-1}\alpha_k)$ such that
			\[
				\Psi_\textup{FB}(z^k+\beta^{-1}\alpha_kd^k)
				-
				\Psi_\textup{FB}(z^k)
				=
				\beta^{-1}\alpha_k\,
					\Psi_\textup{FB}'(z^k+\theta_k d^k)d^k,
			\]
			and together with the above,
			\[
				\Psi_\textup{FB}'(z^k+\theta_kd^k)d^k
				>
				\sigma\,\Psi_\textup{FB}'(z^k)d^k
			\]
			follows. 
			Clearly, $\alpha_k d^k\to 0$ gives $\theta_kd^k\to 0$.
			Noting that $\nabla \Psi_\textup{FB}$ is uniformly continuous on each
			closed ball around $\bar z$, 
			for arbitrary $\varepsilon>0$, we can ensure
			\begin{align*}
				0 
				&<
				\Psi_\textup{FB}'(z^k+\theta_kd^k)d^k
				-
				\sigma\,\Psi_\textup{FB}'(z^k)d^k\\
				&=
				(\Psi_\textup{FB}'(z^k+\theta_kd^k)-\Psi_\textup{FB}'(z^k))d^k
				+
				(1-\sigma)\Psi_\textup{FB}'(z^k)d^k\\
				&\leq
				\varepsilon\nnorm{d^k}
				+
				(1-\sigma)\Psi_\textup{FB}'(z^k)d^k			
			\end{align*}
			for large enough $k\in\N$.
			Combining this with the validity of the angle test gives
			\begin{align*}
				\varepsilon\nnorm{d^k}
				>
				-(1-\sigma)\Psi_\textup{FB}'(z^k)d^k
				\geq
				(1-\sigma)\rho_1\nnorm{\nabla\Psi_\textup{FB}(z^k)}\nnorm{d^k},
			\end{align*}
			i.e., $(1-\sigma)\rho_1\nnorm{\nabla\Psi_\textup{FB}(z^k)}<\varepsilon$
			for all large enough $k\in\N$.
			Since $\varepsilon>0$ has been chosen arbitrarily,
			$\nabla\Psi_\textup{FB}(z^k)\to 0$ follows, which gives
			$\nabla\Psi_\textup{FB}(\bar z)=0$.
		\item Let $\{z^k\}_{k\in K}$ be a subsequence fulfilling $z^k\to_K\bar z$
			for some point $\bar z$ which satisfies the assumptions
			of \cref{lem:local_uniform_invertibility}.
			Furthermore, we note that $\calF_{\textup{max}}$ is locally Lipschitzian by construction.
			Clearly, by $\nnorm{\calF_\textup{FB}(z^k)}\to_K0$, which holds since
			$\bar z$ solves $\calF_\textup{max}(z)=0$ and, thus, also $\calF_\textup{FB}(z)=0$,
			$\{\nu_k\}_{k\in K}$ is a null sequence. Due to
			\cref{lem:equivalence_of_residuals}, it holds
			$\nu_k\in\OO(\nnorm{\calF_{\textup{max}}(z^k)})$,
			so we know that, for sufficiently large $k\in K$, $z^k$ lies in the
			radius of attraction of $\bar z$ mentioned in \cref{thm:local_convergence_LM}
			while $\nu_k$ is sufficiently small in order to apply
			\cref{thm:local_convergence_LM}
			to get the desired results if the LM direction is
			actually accepted.
			This, however, follows for all large enough $k\in K$
			from \cref{cor:quotients_of_residuals} since
			\cref{lem:equivalence_of_residuals} gives
			\[
				\frac{\Psi_\textup{FB}(z^k+d^k)}{\Psi_\textup{FB}(z^k)}
				=
				\left(\frac{\nnorm{\calF_\textup{FB}(z^k+d^k)}}{\nnorm{\calF_\textup{FB}(z^k)}}\right)^2
				\leq
				\left(\frac{K_2}{K_1}\right)^2
				\left(\frac{\nnorm{\calF_{\textup{max}}(z^k+d^k)}}{\nnorm{\calF_{\textup{max}}(z^k)}}\right)^2
				\to_K
				0
			\]
			for the LM directions $\{d^k\}_{k\in K}$.
	\end{enumerate}
\end{proof}

\begin{remark}\label{rem:globalized_LM}\leavevmode
	\begin{enumerate}
		\item Clearly, \cref{alg:GlobalLM} is a descent method with respect to
			$\Psi_\textup{FB}$, i.e., the sequence $\{\nnorm{\calF_\textup{FB}(z^k)}\}_{k\in\N}$
			is monotonically decreasing. This directly gives that $\{\nu_k\}_{k\in\N}$ is
			monotonically decreasing and, thus, bounded by its trivial lower boundedness.
		\item Besides the standard angle test, there is another condition in
			\cref{step:angle_test_+} which avoids that the LM
			direction is chosen if it tends to zero while the angle test is passed.
			This is due to the following observation.
			Suppose that (along a suitable subsequence without relabeling), the
			LM directions $\{d^k\}_{k\in\N}$ pass the angle test but
			tend to zero. In order to prove in \cref{thm:global_convergence_LM} that
			the accumulation points of $\{z^k\}_{k\in\N}$ are stationary for $\Psi_\textup{FB}$,
			one can exploit boundedness of the matrices
			$D_N\calF_{\textup{max}}(z^k)^\top D_N\calF_{\textup{max}}(z^k)+\nu^k\mathbb I^p$ which would give
			\begin{equation}\label{eq:stationarity_of_residual_via_nasty_system}
				\widehat{\mathcal M}_I(\bar w,\bar\xi)
				^\top
				\begin{bmatrix}
					\calH(\bar w,\bar\xi)
					\\
					\max(\calG(\bar w,\bar\xi),-\bar\xi)
				\end{bmatrix}
				=
				0
			\end{equation}
			by definition of the LM direction
			where $(\bar w,\bar\xi)$ is an accumulation point of $\{z^k\}_{k\in\N}$.
			Above, we used the matrix
			\[
				\widehat{\mathcal M}_I(\bar w,\bar\xi)
				:=
				\begin{bmatrix}
					\calH'_w(\bar w,\bar\xi)	&	\calH'_\xi(\bar w,\bar\xi)
					\\
					\mathbb D(I^>(\bar w,\bar\xi)\cup I)\calG'_w(\bar w,\bar\xi)
					&
					\mathbb D(I^>(\bar w,\bar\xi)\cup I)\calG'_\xi(\bar w,\bar\xi)-\mathbb D(I^<(\bar w,\bar\xi)\cup I_\textup{c})
				\end{bmatrix}
			\]	
			as well as the index set
			\[
				I^>(\bar w,\bar \xi):=\{i\in\{1,\ldots,p_2\}\,|\,\calG_i(\bar w,\bar\xi)>-\bar\xi\},
			\]
			and the pair $(I,I_\textup{c})$ is a disjoint partition
			of $\{i\in\{1,\ldots,m\}\,|\,\calG_i(\bar w,\bar \xi)=-\bar\xi_i\}$.
			We note, however, that \eqref{eq:stationarity_of_residual_via_nasty_system} is underdetermined, so we cannot
			deduce $\calH(\bar w,\bar\xi)=0$ and $\max(\calG(\bar w,\bar\xi),-\bar\xi)=0$
			which would give us stationarity of $\bar z:=(\bar w,\bar\xi)$
			for $\Psi_\textup{FB}$.
			This is pretty much in contrast to the situation in
			\cite[Theorem~4.6]{DeLucaFacchineiKanzow2000} where, for square systems,
			stationarity has been shown under some additional assumptions.
		\item Following the literature,
			see e.g.\ \cite{DeLucaFacchineiKanzow2000,FacchineiKanzow1997},
			it is also possible to incorporate inexact solutions of the
			LM equation in \cref{alg:GlobalLM} in a canonical
			way. Combined with a suitable solver for this equation, this approach
			may lead to an immense speed-up of the method. For brevity of presentation,
			we omit this discussion here but just point out the possibility of investigating it.
	\end{enumerate}
\end{remark}

\subsection{On using the Fischer--Burmeister function exclusively}\label{sec:global_LM_FB}

In this subsection, we briefly comment on an algorithm, related to
\cref{alg:GlobalLM}, which fully relies on the residual $\calF_\textup{FB}$ introduced
in \eqref{eq:special_residual_complementarity}.
For completeness, we first present a result regarding the Newton-differentiability of
this function which basically follows from the chain rule stated in \cref{lem:chain_rule_Newton_derivative} 
and \cref{ex:Newton_differentiability}\,\ref{item:newton_differentiability_FB}.

\begin{lemma}\label{lem:Newton_differentiability_of_FB}
	On each bounded subset of $\R^{p_1}\times\R^{p_2}$, the mapping $\calF_\textup{FB}$ is
	Newton-differentiable with Newton-derivative given by
	\[
		(w,\xi)
		\mapsto
		\begin{bmatrix}
			\calH'_w(w,\xi)&\calH'_\xi(w,\xi)\\
			\widetilde{\mathbb D}_\calG(w,\xi)\calG'_w(w,\xi)
				&\widetilde{\mathbb D}_\calG(w,\xi)\calG'_\xi(w,\xi)-\widetilde{\mathbb D}_\xi(w,\xi)
		\end{bmatrix},
	\]
	and this mapping, again, is bounded on bounded sets.
	Whenever the derivatives of $\calG$ and $\calH$ are locally Lipschitzian, 
	the order of Newton-differentiability is $1$.
	Above, the matrices
	$\widetilde{\mathbb D}_\calG(w,\xi),\widetilde{\mathbb D}_\xi(w,\xi)\in\R^{p_2\times p_2}$
	are those ones defined in \cref{lem:gradient_psi_FB}.
\end{lemma}

Subsequently, we will denote the Newton-derivative of $\calF_\textup{FB}$ characterized above
by $D_N\calF_\textup{FB}$. Observe that, due to \cref{lem:gradient_psi_FB}, we have
$D_N\calF_\textup{FB}(z)^\top \calF_\textup{FB}(z)=\nabla\Psi_\textup{FB}(z)$ for each
$z\in\R^p$.

We also need to figure out, in which situations the Newton-derivative from
\cref{lem:Newton_differentiability_of_FB} satisfies the assumptions for local
fast convergence.
\begin{lemma}\label{lem:local_uniform_invertibility_FB}
	Fix a solution $(\bar w,\bar\xi)\in\R^{p_1}\times\R^{p_2}$ of \eqref{eq:NLCS} and
	assume that for each pair $(a,b)\in\R^{p_2}\times\R^{p_2}$ of vectors satisfying
	\begin{equation}\label{eq:ab_vectors_FB_approach}
		\begin{aligned}
			&\forall i\in I^-(\bar w,\bar\xi)\colon&\quad
				&a_i=0,\,b_i=1,\\
			&\forall i\in I^+(\bar w,\bar\xi)\colon&\quad
				&a_i=1,\,b_i=0,\\
			&\forall i\in I^{00}(\bar w,\bar\xi)\colon&\quad
				&(a_i-1)^2+(b_i-1)^2=1,
		\end{aligned}
	\end{equation}
	the matrix
	\begin{equation}\label{eq:matrix__to_check_FB}
		\begin{bmatrix}
			\calH'_w(\bar w,\bar\xi)	&	\calH'_\xi(\bar w,\bar \xi)\\
			\widehat{\mathbb D}^a_\calG(\bar w,\bar\xi)\calG'_w(\bar w,\bar\xi)
			&
			\widehat{\mathbb D}^a_\calG(\bar w,\bar\xi)\calG'_\xi(\bar w,\bar\xi)
				-\widehat{\mathbb D}^b_\xi(\bar w,\bar\xi)
		\end{bmatrix},
	\end{equation}
	where we used
	\begin{align*}
		\widehat{\mathbb D}^a_\calG(\bar w,\bar \xi)&:=\diag(a)\mathbb D(I^+(\bar w,\bar\xi)\cup I^{00}(\bar w,\bar\xi)),\\
		\widehat{\mathbb D}^b_\xi(\bar w,\bar\xi)&:=\diag(b)\mathbb D(I^-(\bar w,\bar\xi)\cup I^{00}(\bar w,\bar\xi)),
	\end{align*}		
	possesses linearly independent columns.
	Then the matrices $\mathcal N(w,\xi)^\top\mathcal N(w,\xi)$, where
	\[
		\mathcal N(w,\xi)
		:=
		\begin{bmatrix}
			\calH'_w(w,\xi)&\calH'_\xi(w,\xi)\\
			\widetilde{\mathbb D}_\calG(w,\xi)\calG'_w(w,\xi)
				&\widetilde{\mathbb D}_\calG(w,\xi)\calG'_\xi(w,\xi)-\widetilde{\mathbb D}_\xi(w,\xi)
		\end{bmatrix},
	\]
	are uniformly positive definite in a neighborhood of $(\bar w,\bar\xi)$.
\end{lemma}
\begin{proof}
	For the proof, we partially mimic our arguments from the proof of
	\cref{lem:local_uniform_invertibility}.
	Thus, let us suppose that the assertion is false.
	Then there exist sequences $\{w^k\}_{k\in\N}\subset\R^{p_1}$,
	$\{\xi^k\}_{k\in\N}\subset\R^{p_2}$,
	$\{(d^k_w,d^k_\xi)\}_{k\in\N}\subset\R^{p_1}\times\R^{p_2}$, and
	$\{\eta_k\}_{k\in\N}\subset[0,\infty)$
	such that $w^k\to\bar w$, $\xi^k\to\bar\xi$, $\eta_k\to 0$, and,
	for each $k\in\N$, $\nnorm{d^k_w}+\nnorm{d^k_\xi}=1$ as well as
	\[
		\eta_k
		=
		\begin{bmatrix}d^k_w\\d^k_\xi\end{bmatrix}^\top
		\mathcal N(w^k,\xi^k)^\top
		\mathcal N(w^k,\xi^k)
		\begin{bmatrix}d^k_w\\d^k_\xi\end{bmatrix}
		.
	\]
	Again, we make use of the abbreviations from \eqref{eq:some_abbrevs}
	and obtain, by definition of the matrix $\mathcal N(w^k,\xi^k)$,
	\begin{equation}\label{eq:surrogate_system_uniform_invertibility_FB}
		\begin{aligned}
		\eta_k
		=
		\nnorm{\calH'_w(k)d^k_w}^2
		&
		+
		2(\calH'_w(k)d^k_w)^\top(\calH'_\xi(k)d^k_\xi)
		+
		\nnorm{\calH'_\xi(k)d^k_\xi}^2
		\\
		&
		+
		\sum_{i=1}^{p_2}(v_a(w^k,\xi^k))_i^2\bigl((\calG_i)'_w(k)d^k_w+(\calG_i)'_\xi(k)d^k_\xi\bigr)^2
		\\
		&
		-
		2\sum_{i=1}^{p_2}
			(v_a(w^k,\xi^k))_i(v_b(w^k,\xi^k))_i\bigl((\calG_i)'_w(k)d^k_w+(\calG_i)'_\xi(k)d^k_\xi\bigr)(d^k_\xi)_i
		\\
		&
		+
		\sum_{i=1}^{p_2}(v_b(w^k,\xi^k))_i^2(d^k_\xi)^2_i
		\end{aligned}
	\end{equation}
	where we exploited the vectors $v_a(w^k,\xi^k)$ and $v_b(w^k,\xi^k)$
	defined in \eqref{eq:surrogate_vectors_FB}.
	For each $i\in I^+(\bar w,\bar\xi)$, $i\in I^{00}_\textup{c}(w^k,\xi^k):=\{1,\ldots,p_2\}\setminus I^{00}(w^k,\xi^k)$ holds for large
	enough $k\in\N$ by continuity of $\mathcal G_i$, and we find the convergences
	$(v_a(w^k,\xi^k))_i\to 1$ and $(v_b(w^k,\xi^k))_i\to 0$ as $k\to\infty$.
	Similarly, for each $i\in I^-(\bar w,\bar\xi)$, we find $i\in I^{00}_\textup{c}(w^k,\xi^k)$ for
	all large enough $k\in\N$, and we also have
	$(v_a(w^k,\xi^k))_i\to 0$ and $(v_b(w^k,\xi^k))_i\to 1$ in this case.
	It remains to consider the indices $i\in I^{00}(\bar w,\bar\xi)$.
	By construction, we know that the sequence 
	$\{((v_a(w^k,\xi^k))_i,(v_b(w^k,\xi^k))_i)\}_{k\in\N}\subset\R^2$
	belongs to the sphere of radius $1$ around $(1,1)$ for each $k\in\N$
	and, thus, possesses an accumulation point in this sphere.
	Thus, taking into account that $I^{00}(\bar w,\bar\xi)$ is a finite set,
	we find an infinite set $K\subset\N$ and vectors $a,b\in\R^{p_2}$
	satisfying $v_a(w^k,\xi^k)\to_K a$, $v_b(w^k,\xi^k)\to_K b$, and
	\eqref{eq:ab_vectors_FB_approach}.
	We may also assume for simplicity that the convergences $d^k_w\to_K d_w$ and $d^k_\xi\to d_\xi$
	hold for a pair $(d_w,d_\xi)\in\R^{p_1}\times\R^{p_2}$ which is non-vanishing.
	Taking the limit $k\to_K\infty$ in 
	\eqref{eq:surrogate_system_uniform_invertibility_FB} then gives
	\begin{align*}
		0
		=
		\nnorm{\calH'_w(\bar w,\bar\xi)d_w}^2
		&
		+
		2(\calH'_w(\bar w,\bar\xi)d_w)^\top(\calH'_\xi(\bar w,\bar\xi)d_\xi)
		+
		\nnorm{\calH'_\xi(\bar w,\bar\xi)d_\xi}^2
		\\
		&+
		\sum_{i\in I^+(\bar w,\bar\xi)}\bigl((\calG_i)'_w(\bar w,\bar\xi)d_w+(\calG_i)'_\xi(\bar w,\bar\xi)d_\xi\bigr)^2
		+
		\sum_{i\in I^-(\bar w,\bar\xi)}(d_\xi)^2_i
		\\
		&+
		\sum_{i\in I^{00}(\bar w,\bar\xi)}
			\bigl(a_i\bigl((\calG_i)'_w(\bar w,\bar\xi)d_w+(\calG_i)'_\xi(\bar w,\bar\xi)d_\xi\bigr)-b_i(d_\xi)_i)^2
	\end{align*}
	which, similar as in the proof of \cref{lem:local_uniform_invertibility}, 
	implies that $(d_w,d_\xi)$ belongs to the kernel
	of the matrix from \eqref{eq:matrix__to_check_FB}.
	This, however, contradicts the lemma's assumptions.
\end{proof}

We note that the assumption of \cref{lem:local_uniform_invertibility_FB}, which corresponds to the full row rank
of all matrices in Bouligand's generalized Jacobian of $\calF_\textup{FB}$ at the reference point, i.e.,
BD-regularity, is
more restrictive than the one from \cref{lem:local_uniform_invertibility}.
Indeed, if the assumption of \cref{lem:local_uniform_invertibility_FB} holds,
then one can choose the vectors $a,b\in\R^{p_2}$ such that
\[
	\forall i\in\{1,\ldots,p_2\}\colon\quad
	a_i:=
	\begin{cases}
		0	& i\in I^-(\bar w,\bar\xi)\cup I_\textup{c},\\
		1	& i\in I^+(\bar w,\bar\xi)\cup I,
	\end{cases}
	\quad
	b_i:=
	\begin{cases}
		0	& i\in I^+(\bar w,\bar\xi)\cup I,\\
		1	& i\in I^-(\bar w,\bar\xi)\cup I_\textup{c}
	\end{cases}
\]
for arbitrary $I\subset I^{00}(\bar w,\bar\xi)$ 
and $I_\textup{c}:=I^{00}(\bar w,\bar\xi)\setminus I$
in order to validate the assumption of \cref{lem:local_uniform_invertibility}.
Clearly, the assumptions of
\cref{lem:local_uniform_invertibility,lem:local_uniform_invertibility_FB}
coincide whenever the biactive set $I^{00}(\bar w,\bar\xi)$ is empty.
This situation is called strict complementarity in the literature.

The subsequent example shows that the assumptions of \cref{lem:local_uniform_invertibility_FB}
can be strictly stronger than those ones of \cref{lem:local_uniform_invertibility}.
\begin{example}\label{ex:FB_assumptions_vs_max_assumptions}
	Let us consider the mixed linear complementarity system
	\[
		w+\xi=0,\,w\geq 0,\,\xi\geq 0,\,w\xi=0
	\]
	which possesses the uniquely determined solution $(\bar w,\bar\xi):=(0,0)$.
	In the context of this example, the functions $\calG,\calH\colon\R\times\R\to\R$
	are given by
	\[
		\forall (w,\xi)\in\R\times\R\colon\quad
		\calG(w,\xi):=-w,\qquad
		\calH(w,\xi):=w+\xi
		.
	\]
	Clearly, the assumption of \cref{lem:local_uniform_invertibility} holds since the
	matrices
	\[
		\begin{bmatrix}
			1	&	1	\\
			-1	&	0
		\end{bmatrix},
		\qquad
		\begin{bmatrix}
			1	&	1	\\
			0	&	-1	
		\end{bmatrix}
	\]
	possess full column rank $2$.
	However, the matrix
	\[
		\begin{bmatrix}
			1	&	1	\\
			-a	&	-b
		\end{bmatrix}
	\]
	possesses column rank $1$ for $a:=b:=1+\sqrt 2/2$, i.e., the assumptions of
	\cref{lem:local_uniform_invertibility_FB} are violated.
\end{example}

From the viewpoint of semismooth solution methods, it also has been observed in 
\cite[Propositions~2.8 and 2.10, Example~2.1]{DeLucaFacchineiKanzow2000} that
the mixed approach via $\calF_\textup{max}$ needs less restrictive assumptions
than the one via $\calF_\textup{FB}$ in order to yield local fast convergence.

Next, we state the globalized nonsmooth LM method for the numerical
solution of \eqref{eq:NLCS} via exclusive use of $\calF_\textup{FB}$ from
\eqref{eq:special_residual_complementarity} in \cref{alg:GlobalLM_FB}.

\begin{algorithm}[Global nonsmooth Levenberg--Marquardt method for mixed nonlinear complementarity systems via
	Fischer--Burmeister function]\leavevmode
	\label{alg:GlobalLM_FB}
	\setcounter{ALC@unique}{0}
	\begin{algorithmic}[1]
		\REQUIRE starting point $z^0\in\R^{p_1}\times\R^{p_2}$,
			parameters $\kappa\in(0,1)$, $\tau_\textup{abs}>0$, $\sigma,\beta\in(0,1)$,
			$\rho>0$, and $\gamma_1,\gamma_2>0$
		\STATE set $k := 0$
		\WHILE{$\nnorm{\calF_\textup{FB}(z^k)}\geq\tau_\textup{abs}$}
		\STATE set $\nu_k:=\min(\gamma_1,\gamma_2\nnorm{\calF_\textup{FB}(z^k)})$ and compute
			$d^k$ as the uniquely determined
			solution of
			\begin{equation}\label{eq:LM_equation_FB}
				\bigl(D_N\calF_\textup{FB}(z^k)^\top D_N\calF_\textup{FB}(z^k)+\nu_k\mathbb I^p\bigr)d
				=
				-\nabla\Psi_\textup{FB}(z^k)
			\end{equation}
		\IF{$\Psi_\textup{FB}(z^k+d^k)\leq \kappa\Psi_\textup{FB}(z^k)$}
			\label{step:ratio_test_FB}
			\STATE set $z^{k+1}:=z^k+d^k$
		\ELSE
			\IF{$\Psi'_\textup{FB}(z^k)d^k>-\rho\nnorm{\nabla\Psi_\textup{FB}(z^k)}\nnorm{d^k}$}
				\STATE set $d^k:=-\nabla\Psi_\textup{FB}(z^k)$
			\ENDIF
			\STATE set $\alpha_k:=\beta^{i_k}$ where $i_k\in\N$ is the smallest positive integer
				such that
				\[
					\Psi_\textup{FB}(z^k+\beta^{i_k}d^k)
					\leq
					\Psi_\textup{FB}(z^k)
					+
					\beta^{i_k}\sigma\,\Psi_\textup{FB}'(z^k)d^k
				\]
			\STATE set $z^{k+1}:=z^k+\alpha_kd^k$
		\ENDIF
		\STATE set $k:=k+1$
		\ENDWHILE
		\RETURN $z^k$
	\end{algorithmic}
\end{algorithm}

Below, we formulate a convergence result which addresses \cref{alg:GlobalLM_FB}.
\begin{theorem}\label{thm:global_convergence_LM_FB}
	Let $\{z^k\}_{k\in\N}$ be a sequence generated by \cref{alg:GlobalLM_FB}.
	\begin{enumerate}
		\item If $\Psi_\textup{FB}(z^k+d^k)\leq \kappa\Psi_\textup{FB}(z^k)$ holds infinitely many
			times in \cref{step:ratio_test_FB},
			then $\{\Psi_\textup{FB}(z^k)\}_{k\in\N}$ is a null sequence and
			each accumulation point of $\{z^k\}_{k\in\N}$ is a solution of
			\eqref{eq:NLCS}.
		\item Each accumulation point of $\{z^k\}_{k\in\N}$ is a stationary point
			of $\Psi_\textup{FB}$.
		\item If an accumulation point $\bar z$ of $\{z^k\}_{k\in\N}$ satisfies
			the assumptions of \cref{lem:local_uniform_invertibility_FB}, then
			the whole sequence $\{z^k\}_{k\in\N}$ converges to $\bar z$
			superlinearly.
			If the derivatives of $\calG$ and $\calH$ are locally Lipschitz continuous
			functions, then the convergence is even quadratic.
	\end{enumerate}
\end{theorem}
\begin{proof}
	The only major difference to the proof of \cref{thm:global_convergence_LM} addresses
	the second statement. More precisely, we need to show that, if the ratio test in
	\cref{step:ratio_test_FB} is violated along the tail of the sequence, and if
	$d^k\to 0$ along a subsequence (without relabeling) while $z^k\to\bar z$ for some
	$\bar z$, then $\bar z$ is stationary for $\Psi_\textup{FB}$.
	As in the proof of \cref{thm:global_convergence_LM}, this is clear if
	$d^k=-\nabla\Psi_\textup{FB}(z^k)$ holds infinitely many times.
	Thus, without loss of generality, let us assume that $d^k$ is the
	LM direction for all $k\in\N$, i.e., the uniquely
	determined solution of \eqref{eq:LM_equation_FB}.
	Boundedness of $\{\nu_k\}_{k\in\N}$ together with \cref{lem:Newton_differentiability_of_FB}
	gives boundedness of 
	the matrices $\{D_N\calF_\textup{FB}(z^k)^\top D_N\calF_\textup{FB}(z^k)+\nu_k\mathbb I^p\}_{k\in\N}$.
	Thus, $d^k\to 0$ gives $\nabla\Psi_\textup{FB}(\bar z)=0$ by continuity
	of $\nabla\Psi_\textup{FB}$ and definition of $d^k$ in
	\eqref{eq:LM_equation_FB}.
\end{proof}

\section{Applications in optimistic bilevel optimization}\label{sec:optimistic_bpp}

\subsection{Model problem and optimality conditions}\label{sec:obpp_modeling}

We consider the so-called standard optimistic bilevel optimization problem
\begin{equation}\label{initialbilev}\tag{OBPP}
\min\limits_{x,y}\quad F(x,y)\quad\text{s.t.} \quad G(x, y) \leq 0, \quad  y\in S(x),
\end{equation}
where $S\colon\R^n\tto\R^m$ is given by
\begin{equation}\label{eq:solution_map}
	\forall x\in\R^n\colon\quad
	S(x):=\argmin\limits_y\{f(x,y)\,|\,y\in Y(x)\}.
\end{equation}
The terminus ``standard'' has been coined in \cite{Zemkoho2016}.
Above, $F,f\colon\R^n\times\R^m\to\R$ are referred to as
the upper- and lower-level objective function, respectively, and assumed to be
twice continuously differentiable.
Furthermore, $G\colon\R^n\times\R^m\to\R^s$ is the twice continuously differentiable
upper-level constraint function, and the set-valued mapping $Y\colon\R^n\tto\R^m$
is given by
\begin{equation*}%\label{eq:lower_level_constraint_map}
	\forall x\in\R^n\colon\quad
	Y(x):=\{y\in\R^m\,|\,g(x,y)\leq 0\},
\end{equation*}
where the describing lower-level constraint function $g\colon\R^n\times\R^m\to\R^t$ 
is also assumed to be twice continuously differentiable.
The component functions of $g$ will be denoted by $g_1,\ldots,g_t$.

We consider the so-called lower-level value function reformulation
\begin{equation}\label{initialvalfuncform0}\tag{VFR}
	\min\limits_{x,y}\quad F(x,y) \quad\text{s.t.} \quad 
		G(x,y)\leq 0, \quad g(x,y)\leq 0,\quad f(x,y)\leq \varphi(x),
\end{equation}
where $\varphi\colon\R^n\to\overline\R:=\R\cup\{\pm\infty\}$ 
is the lower-level value function given by
\begin{equation}\label{eq:value_function}
	\forall x\in\R^n\colon\quad
	\varphi(x) := \inf\limits_y\{f(x,y)\,|\,y\in Y(x)\}.
\end{equation}
It is well known that \eqref{initialbilev} and \eqref{initialvalfuncform0} possess the same
local and global minimizers.

Let us fix a feasible point $(\bar x,\bar y)\in\R^n\times\R^m$ of \eqref{initialbilev}.
Under mild assumptions, one can show that
the lower-level value function $\varphi$ is locally Lipschitz continuous in a
neighborhood of $\bar x$, and its Clarke subdifferential, see \cite{Clarke1990}, obeys the upper estimate
\begin{equation}\label{eq:upper_estimate_subdifferential}
	\partial^\textup{c}\varphi(\bar x)
	\subset
	\{\nabla_xf(\bar x,\bar y)+g'_x(\bar x,\bar y)^\top\hat\nu\,|\,\hat\nu\in\Lambda(\bar x,\bar y)\},
\end{equation}
where $\Lambda(\bar x,\bar y)$ is the lower-level Lagrange multiplier set which comprises
all vectors $\hat\nu\in\R^t$ such that
\begin{align*}
	\nabla_yf(\bar x,\bar y)+g'_y(\bar x,\bar y)^\top\hat\nu&=0,\\
	\hat\nu\geq 0,\quad g(\bar x,\bar y)\leq 0,\quad\hat\nu^\top g(\bar x,\bar y)&=0.
\end{align*}
Let us now assume that $(\bar x,\bar y)$ is already a local minimizer of \eqref{initialbilev} and, thus,
of \eqref{initialvalfuncform0} as well. Again, under some additional assumptions, $(\bar x,\bar y)$ is
a (nonsmooth) stationary point of \eqref{initialvalfuncform0} (in Clarke's sense). 
Keeping the estimate \eqref{eq:upper_estimate_subdifferential} in mind while noting that, by definition
of the lower-level value function, $f(\bar x,\bar y)=\varphi(\bar x)$ is valid,
this amounts to the existence
of $\mu\in\R^s$, $\nu,\hat\nu\in\R^t$, and $\lambda\in\R$ such that
\begin{subequations}\label{eq:KKT_obpp}
	\begin{align}
	\label{eq:KKT_obpp_xy}
	\nabla F(\bar x,\bar y) + G'(\bar x,\bar y)^\top \mu 
		+  g'(\bar x,\bar y)^\top (\nu - \lambda \hat{\nu}) &=0, \\
	\label{eq:KKT_obpp_ll}
	\nabla_y f(\bar x,\bar y) + g'_y(\bar x,\bar y)^\top \hat{\nu} &= 0, \\
	\label{eq:KKT_obpp_alpha}
	\mu\geq 0, \quad G(\bar x,\bar y)\leq 0, \quad \mu^\top G(\bar x,\bar y)&=0,\\
	\label{eq:KKT_obpp_beta}
	\nu\geq 0, \quad g(\bar x,\bar y)\leq 0, \quad \nu^\top g(\bar x,\bar y)&=0,\\
	\label{eq:KKT_obpp_hatbeta}
	\hat{\nu}\geq 0, \quad g(\bar x,\bar y)\leq 0, \quad 
	\hat{\nu}^\top g(\bar x,\bar y)&=0,\\
	\label{eq:KKT_obpp_lambda}
	\lambda&\geq 0.
	\end{align}
\end{subequations}

In the subsequently stated lemma, we postulate some conditions which ensure that local
minimizers of \eqref{initialbilev} indeed are stationary in the sense that there exist
multipliers which solve the system \eqref{eq:KKT_obpp}.
Related results can be found e.g.\ in
\cite{DempeDuttaMordukhovich2007,DempeZemkoho2011,MordukhovichNamPhan2012}.
As we work with slightly different constraint qualifications, 
we provide a proof for the convenience of the reader.

\begin{lemma}\label{lem:CQ_for_stationarity_in_bilevel_programming}
	Fix a local minimizer $(\bar x,\bar y)\in\R^n\times\R^m$ of \eqref{initialbilev}
	such that $\bar x$ is an interior point of $\dom S$.
	The fulfillment of the following conditions implies that there are multipliers
	$\mu\in\R^s$, $\nu,\hat\nu\in\R^t$, and $\lambda\in\R$
	which solve the system \eqref{eq:KKT_obpp}.
	\begin{enumerate}
		\item\label{item:ass_full_convexity} 
			The functions $f$ and $g_1,\ldots,g_t$ are convex in $(x,y)$.
		\item\label{item:ass_LMFCQ}
			Either the functions $f$ and $g$ are affine in $(x,y)$
			or the lower-level Mangasarian--Fromovitz constraint qualification (LMFCQ) 
			\[
				g'_y(\bar x,\bar y)^\top\hat\nu=0,\,
				\hat\nu\geq 0,\,
				\hat\nu^\top g(\bar x,\bar y)=0
				\quad
				\Longrightarrow
				\quad
				\hat\nu=0
			\]
			holds.
		\item\label{item:ass_calmness} 
			The set-valued mapping $\Phi\colon\R\times\R^s\tto\R^n\times\R^m$ given by
			\[
				\forall(r,u)\in\R\times\R^s\colon\quad
				\Phi(r,u):=
				\{(x,y)\in\gph Y\,|\,f(x,y)-\varphi(x)\leq r,\,G(x,y)\leq u\}
			\]
			is calm at $((0,0),(\bar x,\bar y))$.
	\end{enumerate}
\end{lemma}
\begin{proof}
	The imposed convexity assumptions on the lower-level data functions guarantee 
	that $\varphi$ is convex, see e.g.\ \cite[Lemma~2.1]{TaninoOgawa1984}.
	Given that $\bar x$ is an interior point of $\dom S$, 
	we have $|\varphi(x)|<\infty$ for all $x$ in some neighborhood of $\bar x$,
	and therefore $\varphi$ is Lipschitz continuous around $\bar x$. 
	Thus, \eqref{initialvalfuncform0} is a Lipschitz optimization problem around the point $(\bar x, \bar y)$.
	Hence, it follows from \cite[Theorem~4.1]{HenrionJouraniOutrata2002} and \cite[Theorem~6.12]{RockafellarWetts1998} 
	that the calmness of $\Phi$ at $((0,0),(\bar x,\bar y))$ yields the existence of 
	$\mu\in \R^s$, $\lambda \in \R$, $\vartheta\in\partial^\textup{c}\varphi(\bar x)$,
	and $\upsilon\in N_{\gph Y}(\bar x,\bar y)$ 
	such that condition \eqref{eq:KKT_obpp_alpha} holds together with
	\begin{subequations}\label{eq:KKT_obpp_New}
		\begin{align}
			\label{eq:general_opt_cond_New}
					\nabla F(\bar x, \bar y) + G'(\bar x, \bar y)^\top \mu  
						+ \lambda \left(\nabla f(\bar x, \bar y) - (\vartheta,0)\right) 
						+ \upsilon
					&=0,\\
			\label{eq:opt_val_complementarity}
					\lambda\geq 0,\quad f(\bar x,\bar y)-\varphi(\bar x)\leq 0,\quad 
						\lambda\bigl(f(\bar x, \bar y)-\varphi(\bar x)\bigr)
					&=0.
		\end{align}
	\end{subequations}
	Above, $N_{\gph Y}(\bar x, \bar y)$ denotes for the normal cone, in the sense of convex analysis, 
	to the graph of $Y$ (which is a  convex set under the assumptions made) at the point $(\bar x, \bar y)$. 
	Due to the validity of LMFCQ or the fact that $g_1,\ldots,g_t$ are affine, 
	there exists some $\nu\in \R^t$ satisfying \eqref{eq:KKT_obpp_beta} such that $\upsilon=g'(\bar x, \bar y)^\top \nu$,
	see e.g.\ \cite[Theorem~6.14]{RockafellarWetts1998} for the nonlinear case
	(in the linear case, this is a consequence of the well-known Farkas lemma).
	Furthermore, combining the fulfillment of LMFCQ with the full convexity of the lower-level data functions 
	implies that inclusion \eqref{eq:upper_estimate_subdifferential} holds, see e.g.\ \cite[Theorem~2.1]{TaninoOgawa1984}.
	On the other hand, if $f$ and $g$ are affine, then \eqref{eq:upper_estimate_subdifferential}
	holds due to \cite[Proposition~4.1]{YeWu2008}.
	In both situations, we can find some $\hat{\nu}\in \R^t$ 
	such that \eqref{eq:KKT_obpp_ll}, \eqref{eq:KKT_obpp_hatbeta}, and 
	$\vartheta=\nabla_xf(\bar x,\bar y)+g'_x(\bar x,\bar y)^\top\hat\nu$ are satisfied.	
	Plugging this information into \eqref{eq:general_opt_cond_New} gives
	\begin{align*}
		\nabla_xF(\bar x,\bar y)+G'_x(\bar x,\bar y)^\top\mu+g'_x(\bar x,\bar y)^\top(\nu-\lambda\hat\nu)
		&=0,\\
		\nabla_yF(\bar x,\bar y)+G'_y(\bar x,\bar y)^\top\mu+\lambda\nabla_yf(\bar x,\bar y)+g'_y(\bar x,\bar y)^\top\nu
		&=0.
	\end{align*}
	Now, making use of \eqref{eq:KKT_obpp_ll} yields \eqref{eq:KKT_obpp_xy}.
	Finally, it remains to show that \eqref{eq:opt_val_complementarity} reduces to \eqref{eq:KKT_obpp_lambda}. 
	This, however, is obvious as $f(\bar x,\bar y)=\varphi(\bar x)$ holds due to $\bar y\in S(\bar x)$.
\end{proof}

\begin{remark}\label{rem:KKT_obpp}
\leavevmode
\begin{enumerate}
	\item Note that assumption~\ref{item:ass_full_convexity} in \cref{lem:CQ_for_stationarity_in_bilevel_programming} 
		can be replaced by so-called inner semicontinuity of the lower-level optimal solution mapping $S$ from \eqref{eq:solution_map} 
		at $(\bar x, \bar y)$ as this, together with LMFCQ, 
		still yields validity of the estimate \eqref{eq:upper_estimate_subdifferential}
		although $\varphi$ is likely to be not convex in this situation,
		see e.g.\ \cite[Corollary~5.3, Theorem~6.1]{MordukhovichNamPhan2012}.
	\item We note that assumption~\ref{item:ass_calmness} is potentially weaker than the standard conditions used in the literature
		which combine a so-called partial calmness condition, see \cite{YeZhu1995}, with MFCQ-type conditions with respect to the
		upper- and lower-level inequality constraints.
			On the one hand, we admit that a calmness-type assumption on the constraint including the lower-level value function
			are comparatively restrictive, 
			see e.g.\ \cite{HenrionSurowiec2011,KeYaoYeZhang2022,MehlitzMinchenkoZemkoho2020} for discussions,
			so that \eqref{eq:KKT_obpp} can be seen as a reliable necessary optimality condition for \eqref{initialbilev}
			in selected situations only.
			On the other hand, from a numerical perspective, the computational solution of the system \eqref{eq:KKT_obpp}
			turned out to be surprisingly effective in order to determine minimizers of optimistic bilevel optimization,
			see e.g.\ \cite{FliegeTinZemkoho2021,TinZemkoho2021,ZemkohoZhou2021}, and this observation covers situations
			where the assumptions of \cref{lem:CQ_for_stationarity_in_bilevel_programming} are not necessarily satisfied.
	\item In the case where the functions $f$ and $g$ are affine in $(x,y)$, the associated lower-level value function $\varphi$
		is piecewise affine. Hence, whenever $G$ is affine as well, the set-valued mapping $\Phi$ from assumption~\ref{item:ass_calmness}
		is so-called polyhedral, i.e., its graph is the union of finitely many convex polyhedral sets.
		It is well known that such mappings are inherently calm at all points of their graph, see \cite[Proposition~1]{Robinson1981}.
		Thus, assumptions~\ref{item:ass_full_convexity},~\ref{item:ass_LMFCQ}, and~\ref{item:ass_calmness} are satisfied in this
		situation.
	\item A second standard approach to optimality conditions for bilevel optimization problems is based on its so-called
		Karush--Kuhn--Tucker reformulation, but as shown in \cite{ZemkohoZhou2021}, both approaches are, in general, 
		completely different when comparing the resulting optimality conditions and associated constraint qualifications.
		Optimality conditions which are based on the reformulation \eqref{initialvalfuncform0} merely discriminate from each
		other due to different upper estimates for the subdifferential of the optimal value function and the associated
		constraint qualifications which ensure their validity.
		For a detailed comparison, we refer the interested reader to 
		\cite{DempeDuttaMordukhovich2007,DempeMordukhovichZemkoho2014,DempeZemkoho2011}.
\end{enumerate}
\end{remark}

Let us now transfer the necessary optimality conditions from \eqref{eq:KKT_obpp} into
a mixed nonlinear complementarity system. In order to do it in a reasonable way, we need to comment on
the role of the appearing real number $\lambda$. Therefore, let us mention again that
system \eqref{eq:KKT_obpp} is nothing else but the (nonsmooth) Karush--Kuhn--Tucker system
of \eqref{initialvalfuncform0} where the (Clarke) subdifferential of the 
implicitly known function $\varphi$ has been approximated from above by
initial problem data in terms of \eqref{eq:upper_estimate_subdifferential}.
Having this in mind, there are at least two possible interpretation of the
meaning behind $\lambda$. On the one hand, it may represent the Lagrange multiplier
associated with the constraint $f(x,y)-\varphi(x)\leq 0$ in \eqref{initialvalfuncform0}.
Clearly, in order to incorporate optimality for the lower-level problem into the system
\eqref{eq:KKT_obpp}, the multiplier $\hat\nu$ characterized in 
\eqref{eq:KKT_obpp_ll}, \eqref{eq:KKT_obpp_hatbeta}, has to be meaningful, i.e.,
$\lambda$ has to be positive. Similarly, we can interpret $\lambda$ as a partial penalty
parameter which provides local optimality of $(\bar x,\bar y)$ for
\[
	\min\limits_{x,y}\quad F(x,y)+\lambda(f(x,y)-\varphi(x))
	\quad\text{s.t.}\quad G(x,y)\leq 0,\quad g(x,y)\leq 0
\]
whose (nonsmooth) Karush--Kuhn--Tucker system reduces to \eqref{eq:KKT_obpp} under mild assumptions,
and this is the fundamental idea behind the aforementioned concept of partial calmness from \cite{YeZhu1995}.
We note that, whenever a feasible point $(\bar x,\bar y)\in\R^n\times\R^m$ of \eqref{initialbilev}
is a local minimizer of the above partially penalized problem for some $\lambda$, then this
also holds for each larger value of this parameter.
Similarly, the case $\lambda=0$ would not be reasonable here as this means that lower-level
optimality is not a restriction at all. 
Summing up these considerations, we may work with $\lambda>0$.

Subsequently, we will introduce some potential approaches for the reformulation of
\eqref{eq:KKT_obpp} as a mixed nonlinear complementarity system.

\subsubsection{Parametric approach}\label{sec:parametric_approach}

Following the ideas in \cite{TinZemkoho2021}, we suppose that $\lambda>0$ is not
a variable in the system \eqref{eq:KKT_obpp}, but a given parameter which has to
be chosen before the system is solved. Although this approach is challenging
due to the obvious difficulty of choosing an appropriate value for $\lambda$, it
comes along with some theoretical advantages as we will see below.

For a compact notation, let us introduce the block variables
\[
	w:=\begin{bmatrix}x\\y\end{bmatrix}\in\R^{n+m},\qquad
	\xi:=\begin{bmatrix}\mu\\\nu\\\hat\nu\\\end{bmatrix}\in\R^{s+2t}
\]
as well as, for some fixed $\lambda>0$, Lagrangian-type functions 
$L^\textup{op}_\lambda,\ell^\textup{op}\colon\R^{n+m}\times\R^{s+2t}\to\R$ given by
\begin{align*}
	\forall w\in\R^{n+m},\,\forall \xi\in\R^{s+2t}\colon\quad
	L^\textup{op}_\lambda(w,\xi)
	&:=
	F(x,y)+\mu^\top G(x,y)+(\nu-\lambda\hat\nu)^\top g(x,y),\\
	\ell^\textup{op}(w,\xi)
	&:=
	f(x,y)+\hat\nu ^\top g(x,y).
\end{align*}
Setting
\begin{equation}\label{eq:parametric_approach_bpp}
	\calG(w,\xi)
	:=
	\begin{bmatrix} G(x,y)\\g(x,y)\\g(x,y)\end{bmatrix},
	\qquad
	\calH(w,\xi)
	:=
	\begin{bmatrix} \nabla_w L^\textup{op}_\lambda(w,\xi)\\\nabla_y\ell^\textup{op}(w,\xi)\end{bmatrix},
\end{equation}
the solutions of the associated nonlinear complementarity system \eqref{eq:NLCS} are
precisely the solutions of \eqref{eq:KKT_obpp} with a priori given $\lambda$. 
We note that $\calG$ does not depend on the multipliers  in this setting.

Let us now state some assumptions which guarantee local fast convergence of
\cref{alg:GlobalLM,alg:GlobalLM_FB} when applied to the above setting. Therefore, we first need to fix a point 
$(w,\xi)=((x,y),(\mu,\nu,\hat\nu))\in\R^{n+m}\times\R^{s+2t}$ which satisfies the stationarity
conditions \eqref{eq:KKT_obpp} with given $\lambda>0$. 
For such a point, we introduce the following index sets:
\begin{align*}
	I^0_G(x,y)&:=\{i\in\{1,\ldots,s\}\,|\,G_i(x,y)=0\},\\
	I^-_G(x,y)&:=\{i\in\{1,\ldots,s\}\,|\,G_i(x,y)<0\},\\
	I^+_G(x,y,\mu)&:=\{i\in I^0_G(x,y)\,|\,\mu_i>0\},\\
	I^{00}_G(x,y,\mu)&:=\{i\in I^0_G(x,y)\,|\,\mu_i=0\},\\
	I^0_g(x,y)&:=\{i\in\{1,\ldots,t\}\,|\,g_i(x,y)=0\},\\
	I^-_g(x,y)&:=\{i\in\{1,\ldots,t\}\,|\,g_i(x,y)<0\},\\
	I^+_g(x,y,\nu)&:=\{i\in I^0_g(x,y)\,|\,\nu_i>0\},\\
	I^{00}_g(x,y,\nu)&:=\{i\in I^0_g(x,y)\,|\,\nu_i=0\}.
\end{align*}
Furthermore, we make use of the so-called critical subspace defined by
\[
	\mathcal C(w,\xi)
	:=
	\left\{
		\delta w\in\R^{n+m}\,\middle|\,
		\begin{aligned}
			\nabla G_i(x,y)\delta w&=0&&i\in I^+_G(x,y,\mu)\\
			\nabla g_i(x,y)\delta w&=0&&i\in I^+_g(x,y,\nu)\cup I^+_g(x,y,\hat\nu)
		\end{aligned}
	\right\}.
\]
We say that the lower-level linear independence constraint qualification (LLICQ)
holds at $(x,y)$ whenever the gradients
\[
	\nabla_yg_i(x,y)\quad (i\in I^0_g(x,y))
\]
are linearly independent (note that, at the lower-level stage, only $y$ is a variable). 
Analogously, the bilevel linear independence constraint qualification (BLICQ)
is said to hold at $(x,y)$ whenever the gradients
\[
	\nabla G_i(x,y)\quad(i\in I^0_G(x,y)),\qquad\nabla g_i(x,y)\quad (i\in I^0_g(x,y))
\]
are linearly independent.

The following two theorems are inspired by \cite[Theorem~3.3]{FliegeTinZemkoho2021}.
Our first result provides conditions which guarantee validity of the assumptions of
\cref{lem:local_uniform_invertibility} in the setting \eqref{eq:parametric_approach_bpp}.
These assumptions, thus, give local fast convergence of \cref{alg:GlobalLM}
in the present setting.

\begin{theorem}\label{thm:local_fast_convergence_obpp}
	Let $(w,\xi)=((x,y),(\mu,\nu,\hat\nu))\in\R^{n+m}\times\R^{s+2t}$ be a solution of
	the stationarity system \eqref{eq:KKT_obpp} with given $\lambda>0$. 
	Let LLICQ and BLICQ be satisfied at $(x,y)$. 
	Finally, assume that the second-order condition
	\begin{equation}\label{eq:bilevel_soc}
		\forall\delta w\in\mathcal C(w,\xi)\setminus\{0\}\colon\quad
		\delta w^\top\nabla^2_{ww}L^\textup{op}_\lambda(w,\xi)\delta w>0
	\end{equation}
	holds. 
	Then, in the specific setting modeled in \eqref{eq:parametric_approach_bpp}, 
	the assumptions of \cref{lem:local_uniform_invertibility} are valid.
\end{theorem}
\begin{proof}
	For each index sets $I_G\subset I^{00}_G(x,y,\mu)$,
	$I_g\subset I^{00}_g(x,y,\nu)$, and $\hat I_g\subset I^{00}_g(x,y,\hat\nu)$, 
	we need to show that the system comprising
	\begin{subequations}\label{eq:system_obpp}
		\begin{align}
			\label{eq:system_obpp_z_Lop}
			\nabla^2_{ww}L^\textup{op}_\lambda(w,\xi)\delta w
			+
			G'(x,y)^\top\delta\mu
			+
			g'(x,y)^\top(\delta\nu-\lambda\delta\hat\nu)
			&=0,\\
			\label{eq:system_obpp_z_lop}
			(\nabla_{y}\ell^\textup{op})'_w(w,\xi)\delta w+g'_y(x,y)^\top\delta\hat\nu
			&=0
		\end{align}
	\end{subequations}
	as well as the sign conditions
	\begin{subequations}\label{eq:system_obpp_signs_mixed}
		\begin{align}
			\label{eq:system_obpp_G}
			G_i'(x,y)\delta w
			&=
			0
			\quad
			i\in I^+_G(x,y,\mu)\cup I_G,\\
			\label{eq:system_obpp_g}
			g_i'(x,y)\delta w
			&=
			0
			\quad
			i\in I^+_g(x,y,\nu)\cup I^+_g(x,y,\hat\nu)\cup I_g\cup\hat I_g,\\
			\label{eq:system_obpp_alpha}
			\delta\mu_i
			&=
			0
			\quad
			i\in I^-_G(x,y)\cup(I_G)_\textup{c},\\
			\label{eq:system_obpp_beta}
			\delta\nu_i
			&=
			0
			\quad
			i\in I^-_g(x,y)\cup(I_g)_\textup{c},\\
			\label{eq:system_obpp_hatbeta}
			\delta\hat\nu_i
			&=
			0
			\quad
			i\in I^-_g(x,y)\cup(\hat I_g)_\textup{c}
		\end{align}
	\end{subequations}
	only possesses the trivial solution.
	Above, we used $(I_G)_\textup{c}:=I^{00}_G(x,y,\mu)\setminus I_G$,
	$(I_g)_\textup{c}:=I^{00}_g(x,y,\nu)\setminus I_g$, and
	$(\hat I_g)_\textup{c}:=I^{00}_g(x,y,\hat\nu)\setminus\hat I_g$.
	
	Multiplying \eqref{eq:system_obpp_z_Lop} from the left with $\delta w^\top$ and respecting
	\eqref{eq:system_obpp_signs_mixed} gives
	$\delta w^\top\nabla^2_{ww}L^\textup{op}_\lambda(w,\xi)\delta w=0$.
	On the other hand, \eqref{eq:system_obpp_G} and \eqref{eq:system_obpp_g} yield
	$\delta w\in\mathcal C(w,\xi)$. 
	Hence, the assumptions of the theorem can be used to find $\delta w=0$.
	Now, we can exploit LLICQ in order to obtain $\delta\hat\nu=0$ from \eqref{eq:system_obpp_z_lop}
	and \eqref{eq:system_obpp_hatbeta}.
	Finally, with the aid of BLICQ, \eqref{eq:system_obpp_z_Lop}, \eqref{eq:system_obpp_alpha},
	and \eqref{eq:system_obpp_beta}, we find
	$\delta\mu=0$ and $\delta\nu=0$.
	This shows the claim.
\end{proof}

	Let us note from the proof that \cref{thm:local_fast_convergence_obpp} remains correct whenever 
	\eqref{eq:bilevel_soc} is replaced by the slightly weaker condition
	\begin{equation}\label{eq:bilevel_soc_weak}
		\forall\delta w\in\mathcal C(w,\xi)\setminus\{0\}\colon\quad
		\delta w^\top\nabla^2_{ww}L^\textup{op}_\lambda(w,\xi)\delta w\neq 0.
	\end{equation}
	However, \eqref{eq:bilevel_soc} is related to second-order sufficient optimality
	conditions for the characterization of strict local minimizers of
	\eqref{initialbilev}, see e.g.\ \cite{MehlitzZemkoho2021}, 
	and as we aim to find local minimizers of \eqref{initialbilev}, 
	this seems to be a more natural assumption than \eqref{eq:bilevel_soc_weak}
	which seemingly concerns saddle points of $L^\textup{op}_\lambda$.

Under slightly stronger conditions, we can prove that even the assumptions of
\cref{lem:local_uniform_invertibility_FB} hold in the precise setting from \eqref{eq:parametric_approach_bpp}
which, in turn, guarantee local fast convergence of \cref{alg:GlobalLM_FB} in the present setting.

\begin{theorem}
	\label{thm:local_fast_convergence_obpp_FB}
	Let $(w,\xi)=((x,y),(\mu,\nu,\hat\nu))\in\R^{n+m}\times\R^{s+2t}$ be a solution of
	the stationarity system \eqref{eq:KKT_obpp} with given $\lambda>0$. 
	Let LLICQ and BLICQ be satisfied at $(x,y)$, 
	and let $I^{00}_g(x,y,\hat\nu)\subset I^+_g(x,y,\nu)$ hold.
	Finally, assume that the second-order condition
	\eqref{eq:bilevel_soc} holds. 
	Then, in the specific setting modeled in \eqref{eq:parametric_approach_bpp}, 
	the assumptions of \cref{lem:local_uniform_invertibility_FB} are valid.
\end{theorem}
\begin{proof}
	For each $i\in I^{00}_G(x,y,\mu)$, let $(a_i^G,b_i^G)\in\R^2$ satisfy
	\[
		(a^G_i-1)^2+(b^G_i-1)^2=1.
	\]
	Similarly, for each $i\in I^{00}_g(x,y,\nu)$ ($i\in I^{00}_g(x,y,\hat\nu)$),
	let $(a_i^g,b_i^g)\in\R^2$ ($(\hat a_i^g,\hat b_i^g)\in\R^2$) satisfy
	\[
		(a^g_i-1)^2+(b^g_i-1)^2=1
		\quad
		\bigl((\hat a^g_i-1)^2+(\hat b^g_i-1)^2=1\bigr).
	\]
	We need to show that the system comprising the conditions from
	\eqref{eq:system_obpp}
	as well as the sign conditions
	\begin{subequations}\label{eq:system_obpp_signs}
		\begin{align}
			\label{eq:system_obpp_G_FB}
			G_i'(x,y)\delta w
			&=
			0
			\quad
			i\in I^+_G(x,y,\mu),\\
			\label{eq:system_obpp_g_FB}
			g_i'(x,y)\delta w
			&=
			0
			\quad
			i\in I^+_g(x,y,\nu)\cup I^+_g(x,y,\hat\nu),\\
			\label{eq:system_obpp_G_00}
			a_i^GG'_i(x,y)\delta w-b_i^G\delta\mu_i
			&=
			0
			\quad
			i\in I^{00}_G(x,y,\mu),\\
			\label{eq:system_obpp_g_00}
			a_i^gg'_i(x,y)\delta w-b_i^g\delta\nu_i
			&=
			0
			\quad
			i\in I^{00}_g(x,y,\nu),\\
			\label{eq:system_obpp_g_00_hat}
			\hat a_i^gg'_i(x,y)\delta w-\hat b^g_i\delta\hat\nu_i
			&=
			0
			\quad
			i\in I^{00}_g(x,y,\hat\nu),\\
			\label{eq:system_obpp_alpha_FB}
			\delta\mu_i
			&=
			0
			\quad
			i\in I^-_G(x,y),\\
			\label{eq:system_obpp_beta_FB}
			\delta\nu_i=\delta\hat\nu_i
			&=
			0
			\quad
			i\in I^-_g(x,y)
		\end{align}
	\end{subequations}
	only possesses the trivial solution.
	
	For later use, we introduce index sets $P^G_{10},P^G_{01},P^G_{+}\subset I^{00}_G(x,y,\mu)$ by means of
	\begin{align*}
		P^G_{10}&:=\{i\in I^{00}_G(x,y,\mu)\,|\,a^G_i=1,\,b^G_i=0\},\\
		P^G_{01}&:=\{i\in I^{00}_G(x,y,\mu)\,|\,a^G_i=0,\,b^G_i=1\},\\
		P^G_{+}&:=I^{00}_G(x,y,\mu)\setminus(P^G_{10}\cup P^G_{01}).
	\end{align*}
	Let us note that $a^G_i,b^G_i>0$ holds for each $i\in P^G_+$, which gives
	$G'_i(x,y)\delta w=(b^G_i/a^G_i)\delta\mu_i$ by \eqref{eq:system_obpp_G_00}.
	Furthermore, for each $i\in P^G_{10}$, we have $G'_i(x,y)\delta w=0$,
	while $\delta\mu_i=0$ is valid for each $i\in P^{G}_{01}$.
	Let the index sets $P^g_{10},P^g_{01},P^g_+\subset I^{00}_g(x,y,\nu)$ and
	$\hat P^g_{10},\hat P^g_{01},\hat P^g_+\subset I^{00}_g(x,y,\hat\nu)$ be defined in analogous fashion.
	Note that for each $i\in \hat P^g_+\subset I^{00}_g(x,y,\hat\nu)\subset I^+_g(x,y,\nu)$,
	\eqref{eq:system_obpp_g_FB} and \eqref{eq:system_obpp_g_00_hat} give $(\hat b^g_i/\hat a^g_i)\delta\hat\nu_i = g'_i(x,y)\delta w=0$.
	
	Clearly, \eqref{eq:system_obpp_G_FB} and \eqref{eq:system_obpp_g_FB} give $\delta w\in\mathcal C(w,\xi)$.
	Multiplying \eqref{eq:system_obpp_z_Lop} from the left with $\delta w^\top$ while respecting
	\eqref{eq:bilevel_soc}, \eqref{eq:system_obpp_signs}, and the above discussion gives
	\begin{align*}
		0
		\leq
		\delta w^\top\nabla^2_{ww}L^\textup{op}_\lambda(w,\xi)\delta w
		&=
		-\sum_{i\in P^G_+}\frac{b^G_i}{a^G_i}(\delta\mu_i)^2
		-\sum_{i\in P^g_+}\frac{b^g_i}{a^g_i}(\delta\nu_i)^2
		+\lambda\sum_{i\in\hat P^g_+}\frac{\hat b^g_i}{\hat a^g_i}(\delta\hat\nu_i)^2\\
		&=
		-\sum_{i\in P^G_+}\frac{b^G_i}{a^G_i}(\delta\mu_i)^2
		-\sum_{i\in P^g_+}\frac{b^g_i}{a^g_i}(\delta\nu_i)^2
		\leq 
		0.
	\end{align*}
	Thus, we have $\delta w^\top\nabla^2_{ww}L^\textup{op}_\lambda(w,\xi)\delta w=0$, so \eqref{eq:bilevel_soc}
	yields $\delta w=0$.
	Now, we can argue as in the proof of \cref{thm:local_fast_convergence_obpp} in order to
	obtain $\delta\mu=0$, $\delta\nu=0$, and $\delta\hat\nu=0$ from
	\eqref{eq:system_obpp_z_lop}, \eqref{eq:system_obpp_alpha_FB}, and
	\eqref{eq:system_obpp_beta_FB} with the aid of LLICQ and BLICQ.
	This completes the proof.
\end{proof}

Let us note that both \cref{thm:local_fast_convergence_obpp,thm:local_fast_convergence_obpp_FB} 
drastically enhance \cite[Theorem~3.3]{FliegeTinZemkoho2021} where strict complementarity is assumed.

\subsubsection{Variable approach}\label{sec:variable_approach}

The discussions in \cite{FischerZemkohoZhou2021,TinZemkoho2021} underline that,
	from a numerical point of view,
treating $\lambda$ as a parameter in \eqref{eq:KKT_obpp} is nontrivial since the particular
choice of it is quite involved. We therefore aim to interpret $\lambda$ as a variable which is
determined by the solution algorithm.
In this section, we suggest two associated approaches.

First, we may consider the block variables
\[
	w:=\begin{bmatrix}x\\y\end{bmatrix}\in\R^{n+m},\qquad
	\xi:=\begin{bmatrix}\mu\\\nu\\\hat\nu\\\lambda\end{bmatrix}\in\R^{s+2t+1}
\]
as well as Lagrangian-type functions 
$L^\textup{op}_1,\ell^\textup{op}_1\colon\R^{n+m}\times\R^{s+2t+1}\to\R$ given by
\begin{align*}
	\forall w\in\R^{n+m},\,\forall \xi\in\R^{s+2t+1}\colon\quad
	L^\textup{op}_1(w,\xi)
	&:=
	F(x,y)+\mu^\top G(x,y)+(\nu-\lambda\hat\nu)^\top g(x,y),\\
	\ell^\textup{op}_1(w,\xi)
	&:=
	f(x,y)+\hat\nu^\top g(x,y).
\end{align*}
Setting
\begin{equation}\label{eq:variable_approach_bpp_max}
	\calG(w,\xi)
	:=
	\begin{bmatrix} G(x,y)\\g(x,y)\\g(x,y)\\0\end{bmatrix},
	\qquad
	\calH(w,\xi)
	:=
	\begin{bmatrix} \nabla_w L^\textup{op}_1(w,\xi)\\\nabla_y\ell^\textup{op}_1(w,\xi)\end{bmatrix},
\end{equation}
we can reformulate \eqref{eq:KKT_obpp} as a mixed complementarity system of type \eqref{eq:NLCS}.
Here, $\lambda$ is treated as part of the Lagrange multiplier vector associated with
\eqref{eq:KKT_obpp}, and this seems to be rather natural.

A second idea is based on a squaring trick.
Observe that the system \eqref{eq:KKT_obpp} can be equivalently reformulated by combining
\[
	\nabla F(\bar x,\bar y) + G'(\bar x,\bar y)^\top \mu 
		+  g'(\bar x,\bar y)^\top (\nu - \zeta^2 \hat{\nu}) =0
\]
with \eqref{eq:KKT_obpp_ll}-\eqref{eq:KKT_obpp_hatbeta} for multipliers 
$\mu\in\R^s$, $\nu,\hat\nu\in\R^t$, and $\zeta\in\R$.
This eliminates the sign condition \eqref{eq:KKT_obpp_lambda}.
Thus, using the block variables
\[
	w:=\begin{bmatrix}x\\y\\\zeta\end{bmatrix}\in\R^{n+m+1},\qquad
	\xi:=\begin{bmatrix}\mu\\\nu\\\hat\nu\end{bmatrix}\in\R^{s+2t}
\]
as well as Lagrangian-type functions 
$L^\textup{op}_2,\ell^\textup{op}_2\colon\R^{n+m+1}\times\R^{s+2t}\to\R$ given by
\begin{align*}
	\forall w\in\R^{n+m+1},\,\forall \xi\in\R^{s+2t}\colon\quad
	L^\textup{op}_2(w,\xi)
	&:=
	F(x,y)+\mu^\top G(x,y)+(\nu-\zeta^2\hat\nu)^\top g(x,y),\\
	\ell^\textup{op}_2(w,\xi)
	&:=
	f(x,y)+\hat\nu^\top g(x,y),
\end{align*}
we can recast the resulting system in the form \eqref{eq:NLCS} when using
\begin{equation}\label{eq:variable_approach_bpp_square}
	\calG(w,\xi)
	:=
	\begin{bmatrix} G(x,y)\\g(x,y)\\g(x,y)\end{bmatrix},
	\qquad
	\calH(w,\xi)
	:=
	\begin{bmatrix} \nabla_{(x,y)} L^\textup{op}_2(w,\xi)\\\nabla_y\ell^\textup{op}_2(w,\xi)\end{bmatrix}.
\end{equation}
Note that in both settings \eqref{eq:variable_approach_bpp_max} and \eqref{eq:variable_approach_bpp_square},
the mapping $\calG$ does not depend on the multiplier.

Let us point the reader's attention to the fact that both approaches discussed above come along with the heavy disadvantage
that a result similar to \cref{thm:local_fast_convergence_obpp} does not seem to hold.
In order to see this, one can try to mimic the proof of \cref{thm:local_fast_convergence_obpp}
in the setting \eqref{eq:variable_approach_bpp_max} and \eqref{eq:variable_approach_bpp_square},
respectively. After all possible eliminations, one ends up with $\delta\nu=\delta\lambda\,\hat\nu$
and $\delta\nu=2\zeta\,\delta\zeta\,\hat\nu$, respectively, and both of these linear equations possess a 
nonvanishing kernel. Clearly, these arguments extend to \cref{thm:local_fast_convergence_obpp_FB} as well.

\subsection{Computational experiments}\label{sec:computational_experiments}

We would like to point the reader's attention to the fact that the numerical application of
Gauss--Newton or LM methods for the numerical solution of a smoothed version of system \eqref{eq:KKT_obpp}
has been investigated in \cite{FliegeTinZemkoho2021,TinZemkoho2021}.
Therein, the authors challenged their methods by running experiments on the 124 nonlinear bilevel optimization
problems which are contained in the BOLIB library from \cite{ZhouZemkohoTin2020}.
Let us note that only some of these examples satisfy the assumptions of \cref{lem:CQ_for_stationarity_in_bilevel_programming},
leading to the fact that \eqref{eq:KKT_obpp} does not provide a reliable necessary optimality condition, i.e.,
system \eqref{eq:KKT_obpp} does not possess a solution in many situations.
In this regard, it is not surprising that, in \cite[Section~5.4]{FliegeTinZemkoho2021}, 
the authors admit that the methods under consideration most often
terminate since the maximum number of iterations is reached or since the norm of the difference of two consecutive iterates
becomes small, which is clearly not satisfying.
In \cite[Section~4.1.4]{TinZemkoho2021}, the authors introduce additional termination criteria, based on the difference
of the norm of the FB residual and certain thresholds of the iteration
number, to enforce termination even in situations where the algorithm did not fully converge.
It is not clarified in \cite{TinZemkoho2021} why the authors do not check for approximate stationarity of the squared
FB residual.
Summarizing this, the experiments in \cite{FliegeTinZemkoho2021,TinZemkoho2021} visualize some global properties of the
promoted solution approaches but do not respect the assumptions needed to ensure (local fast) convergence.
The results are, thus, of limited expressiveness.

Here, we restrict ourselves to the illustration of certain features of our nonsmooth LM methods from
\cref{alg:GlobalLM,alg:GlobalLM_FB} on problems where the assumptions of \cref{lem:CQ_for_stationarity_in_bilevel_programming} 
or even \cref{thm:local_fast_convergence_obpp,thm:local_fast_convergence_obpp_FB} hold.
Furthermore, we take care of distinguishing between the parametric approach from \cref{sec:parametric_approach},
where the multiplier $\lambda$ is treated as a parameter which has to be chosen a priori, 
see \cite{FischerZemkohoZhou2021,FliegeTinZemkoho2021} for further comments on how precisely this
parameter can be chosen in numerical practice,
and the variable approach from \cref{sec:variable_approach}, where $\lambda$ is an additional variable.

\subsubsection{Implementation and setting}

We implemented \cref{alg:GlobalLM} (method \textbf{mixLM}) and \cref{alg:GlobalLM_FB} (method \textbf{FBLM}) 
for the three different settings described in
\eqref{eq:parametric_approach_bpp} (setting \textbf{Para}), 
\eqref{eq:variable_approach_bpp_max} (setting \textbf{Var1}), and 
\eqref{eq:variable_approach_bpp_square} (setting \textbf{Var2})
(making a total number of six algorithms) in MATLAB2022b based on user-supplied gradients and Hessians.  
Numerical experiments were ran for two bilevel optimization problems
whose minimizers satisfy the stationarity conditions from \eqref{eq:KKT_obpp}:
\begin{itemize}[leftmargin=9em]
	\item[\textbf{Experiment~1}] the problem from \cite[Example~8]{MehlitzZemkoho2021} where \eqref{eq:KKT_obpp}
		holds at the global minimizer for each $\lambda>0$,
	\item[\textbf{Experiment~2}] the inverse transportation problem from \cite[Section~5.3.2, Experiment~3]{Mehlitz2022}
		whose lower-level problem and upper-level constraints are fully affine.
\end{itemize}

For each experiment, a certain number of (random) starting points has been generated, 
and each of the algorithms has been run based on these starting points.
If the termination criterion $\nnorm{\mathcal F_\textup{FB}(z^k)}<\tau_\textup{abs}$ (\textbf{Term$=$1}) is violated in
\cref{alg:GlobalLM,alg:GlobalLM_FB}, we additionally check validity of $\nnorm{\nabla\Psi_\textup{FB}(z^k)}<\tau_\textup{abs}^\textup{stat}$
(\textbf{Term$=$2}) for some constant $\tau_\textup{abs}^\textup{stat}>0$ in order to capture situations 
where a stationary point of the squared FB residual $\Psi_\textup{FB}$ has been found.
Furthermore, we equip each algorithm with a maximum number of possible iterations, and terminate if it is hit (\textbf{Term$=$0}).

In order to compare the output of the six methods, we make use of so-called performance profiles, see \cite{DolanMore2002},
based on the following indicators: total number of iterations, execution time (in seconds), upper-level objective value, and 
percentage of full LM steps (i.e., the quotient \# full LM steps/\# total iterations).
We denote by $t_{s,i}>0$ the metric of comparison (associated with the current performance index) 
of solver $s \in \mathcal S$ for solving the instance $i \in \mathcal I$ of the problem, 
where $\mathcal S$ is the set of solvers and $\mathcal I$ represents the different starting points. 
We define the performance ratio by 
\begin{equation*}
	\forall s\in\mathcal S,\,\forall i\in\mathcal I\colon\quad
    r_{s,i} := \dfrac{t_{s,i}}{\min\{t_{s',i}\,|\, s' \in \mathcal S\}}.
\end{equation*}
This quantity is the ratio of the performance of solver $s\in \mathcal S$ to solve instance $i\in \mathcal I$ 
compared to the best performance of any other algorithm in $\mathcal S$ to solve instance $i$. 
The cumulative distribution function $\omega_s\colon[1,\infty)\to[0,1]$ of the current performance index
associated with the solver $s\in\mathcal S$ is defined by
\begin{equation*}
    \forall\tau\in[1,\infty)\colon\quad
    \omega_{s}(\tau) := \frac{|\{i \in \mathcal I\,|\, r_{s,i} \leq \tau\}|}{|\mathcal I|}.
\end{equation*}
The performance profile for a fixed performance index shows
(the illustrative parts of) the graphs of all the functions $\omega_{s}$, $s\in \mathcal S$.
The value $\omega_{s}(1)$ represents the fraction of problem instances for which solver $s\in\mathcal S$ shows the best performance. 
For arbitrary $\tau\geq 1$, $\omega_s(\tau)$ shows the fraction of problem instances for which solver $s\in\mathcal S$ 
shows at most the $\tau$-fold of the best performance.

Finally, let us comment on the precise construction of the performance metric. 
For the comparison of the total number of iterations and execution time (in seconds), we simply rely on the index under consideration.
Computed upper-level objective values are shifted by the minimum function value known in the literature to ensure that the
associated performance metric does not produce negative outputs 
(additionally, we also add a small positive offset to ensure positivity). 
In order to benchmark the percentage of full LM steps, we subtract the aforementioned quotient of
full LM steps over total number of steps from $1$ (again, a positive offset is added to ensure positivity).
This guarantees that small values of the metric indicate a desirable behavior.

\subsubsection{Numerical examples}

\paragraph{Experiment~1}

We investigate the bilevel optimization problem from \cite[Example~8]{MehlitzZemkoho2021} which is given by
\[
	\min\limits_{x,y}\quad (x-8)^2+(y-9)^2
	\quad\text{s.t.}\quad 
	x\geq 0, \quad
	y\in\argmin\limits_y\{(y-3)^2\,|\,y^2\leq x\}.
\]
Its global minimizer is $(\bar x,\bar y):=(9,3)$, and the corresponding stationarity conditions from 
\eqref{eq:KKT_obpp} hold for each $\lambda>0$ when choosing $\mu:=0$, $\hat\nu\geq 0$, and
$\nu:=2+\lambda\hat\nu$.
One can easily check that the assumptions of \cref{thm:local_fast_convergence_obpp} are valid.
If $\hat\nu>0$ is chosen, strict complementarity holds.
However, even for $\hat\nu:=0$, we have $I^{00}_g(\bar x,\bar y,\hat\nu)=I^+_g(\bar x,\bar y,\nu)$
so that the assumptions of \cref{thm:local_fast_convergence_obpp_FB} hold as well for each feasible
choice of the multipliers.

\begin{figure}[ht]
	\begin{center}
		\includegraphics[height=5.0cm]{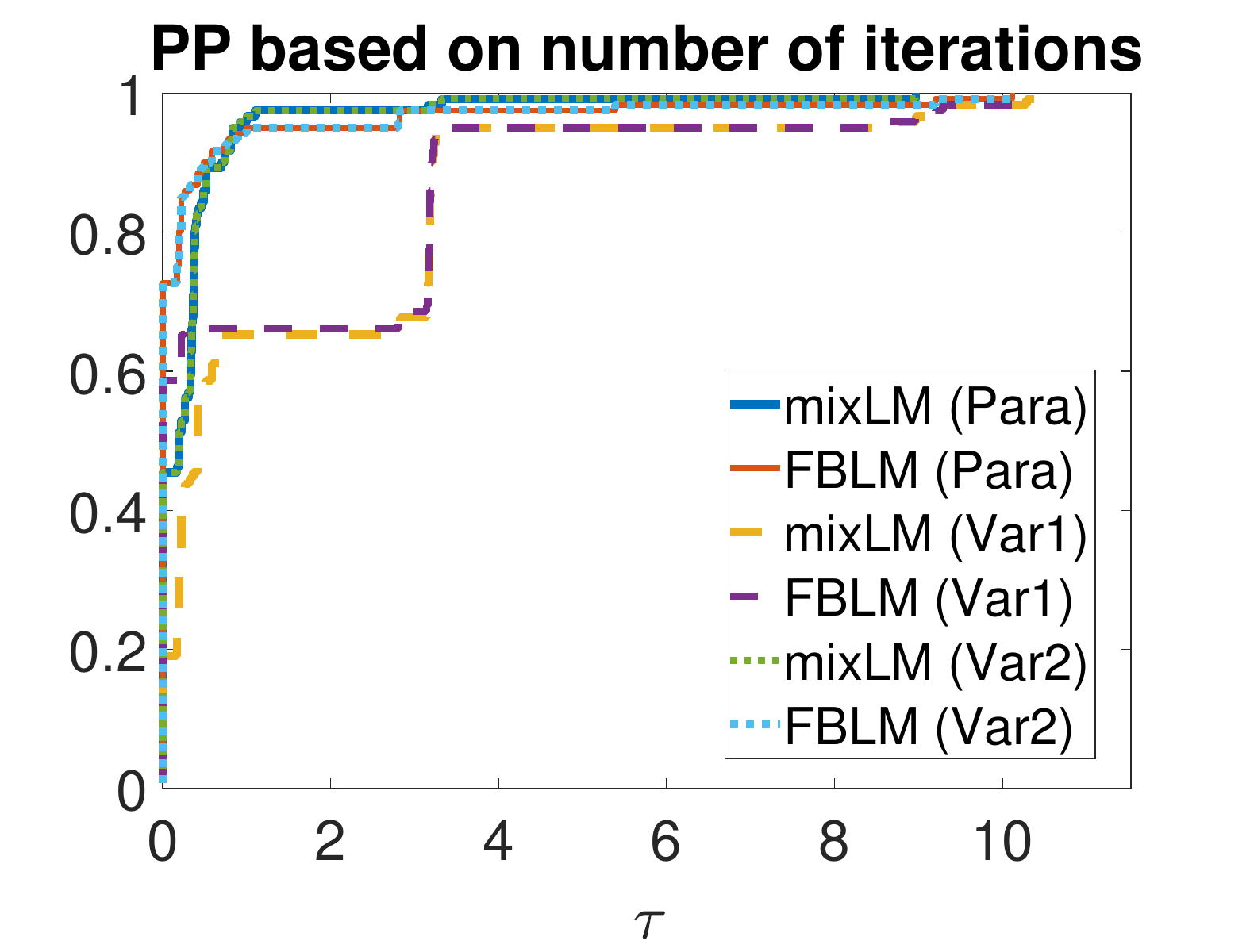} 
		\includegraphics[height=5.0cm]{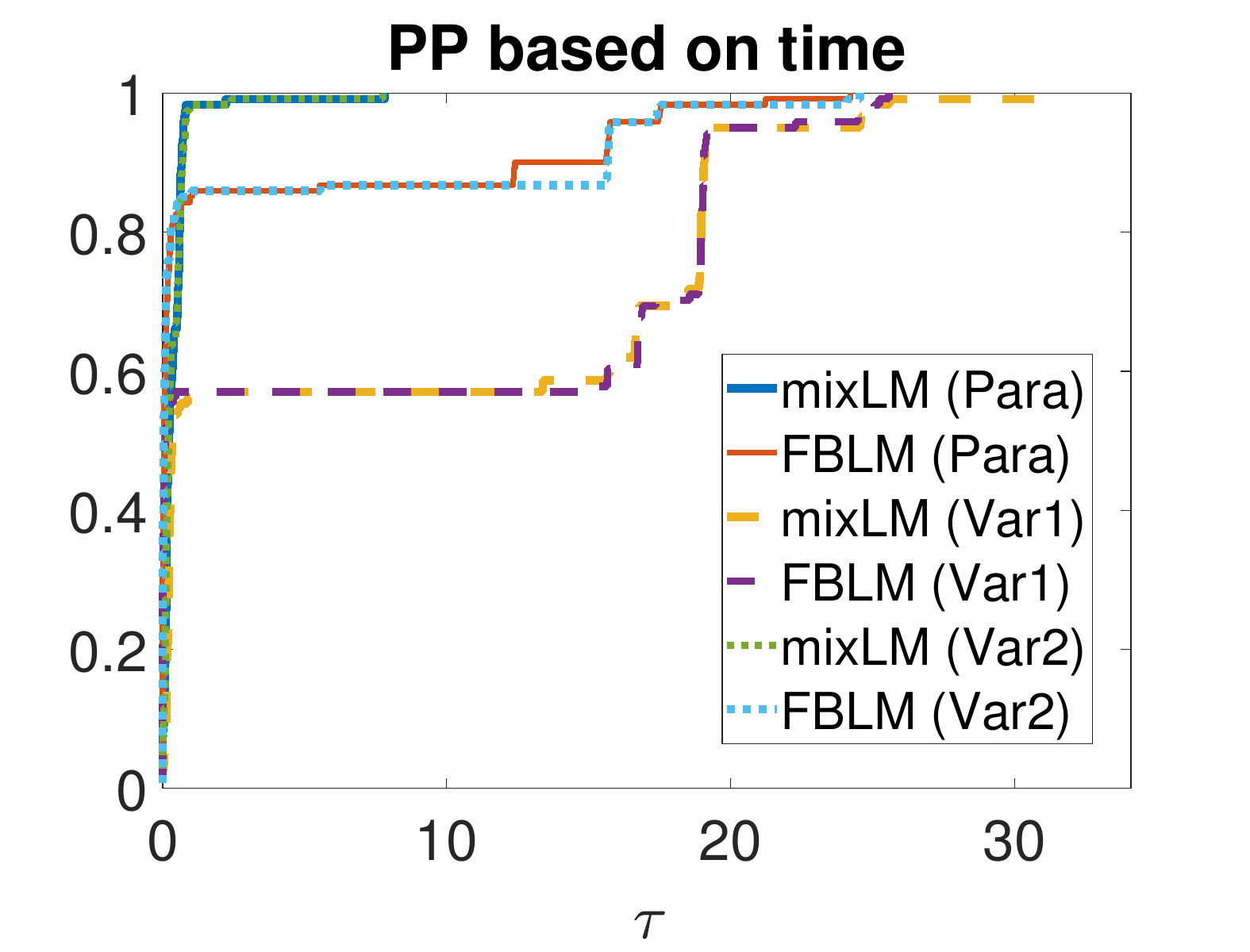}\\ 
		\includegraphics[height=5.0cm]{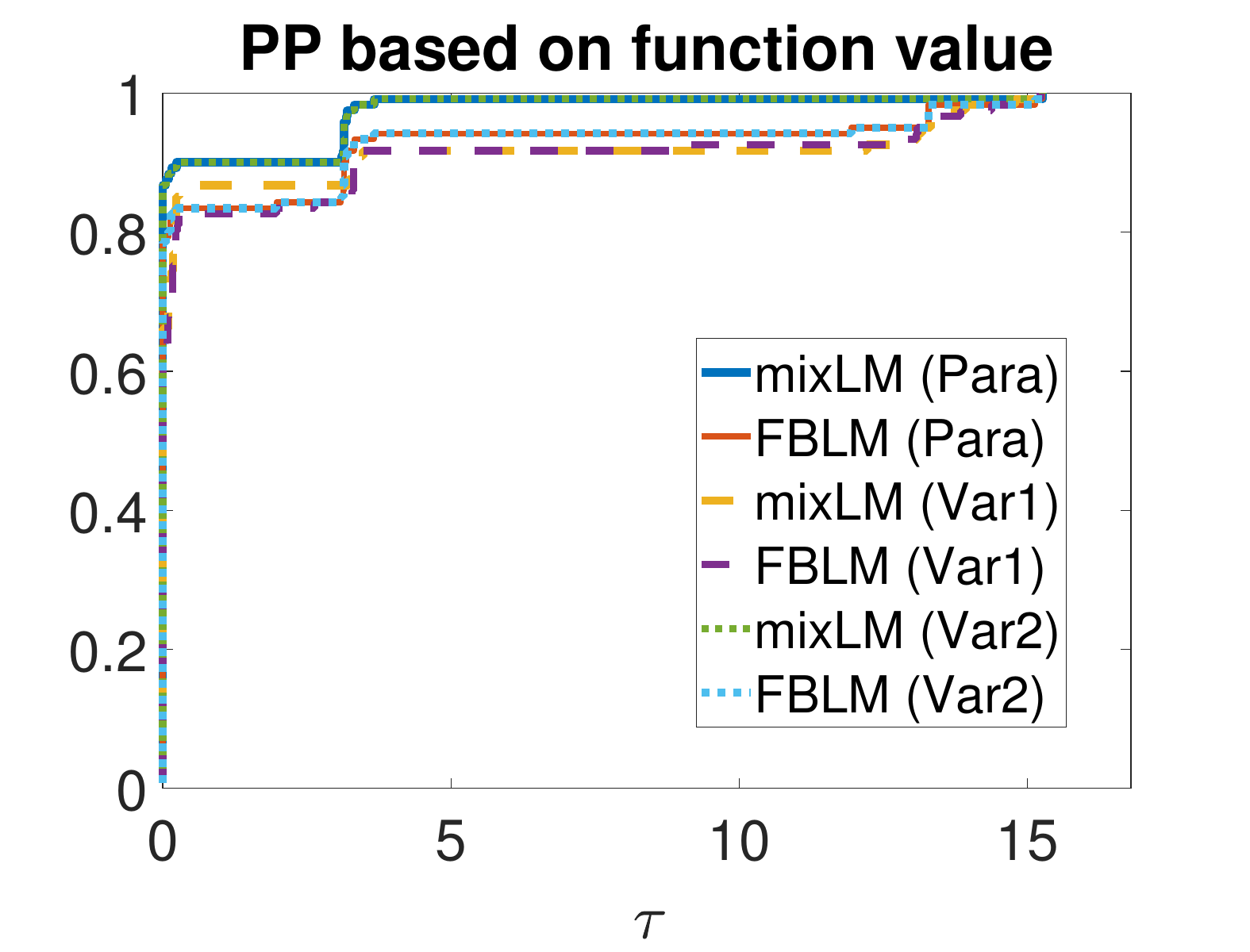} 
		\includegraphics[height=5.0cm]{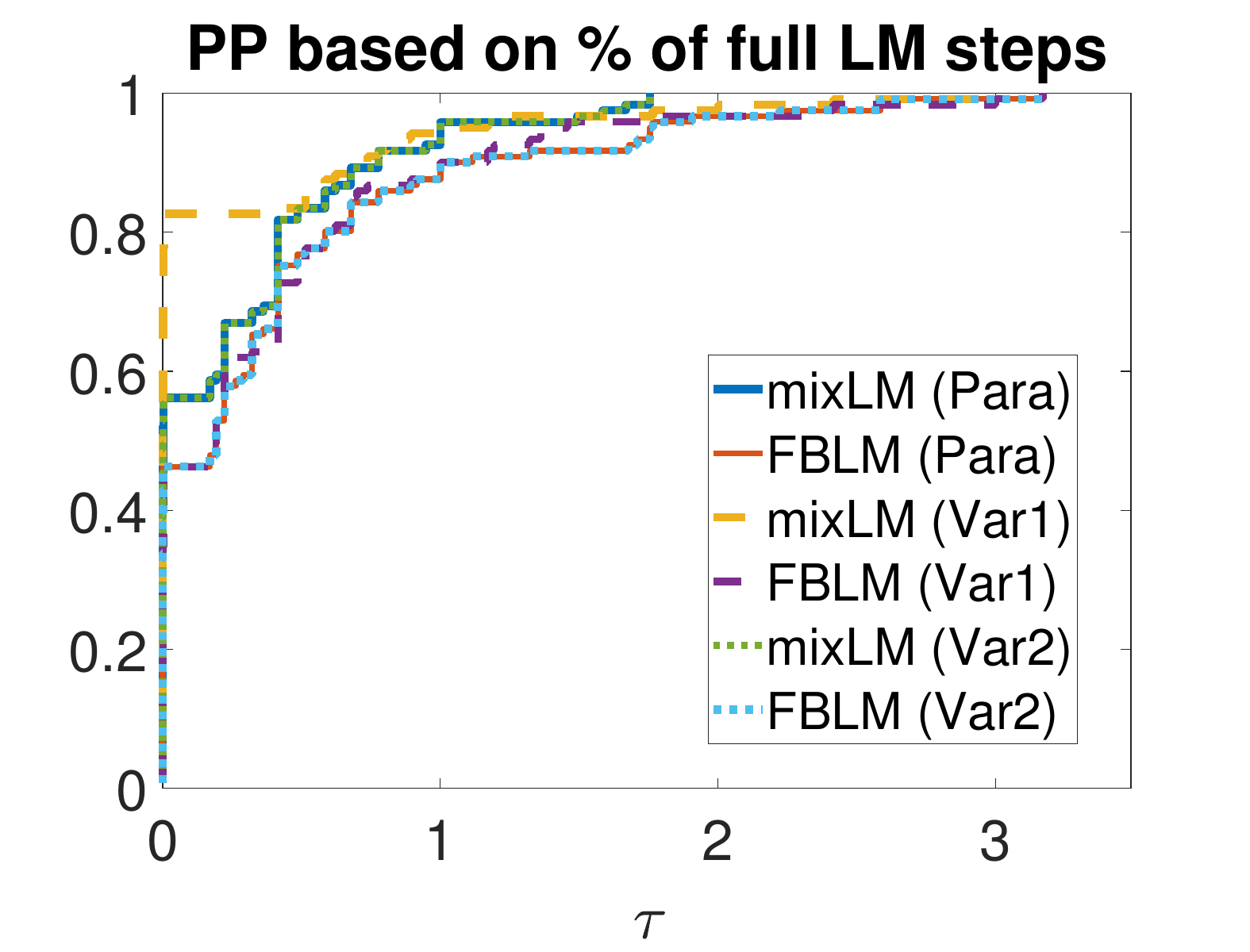} 
	\end{center}
	\caption{%
		Performance profiles for Experiment~1.
		From top left to bottom right: 
    	number of total iterations, execution time, upper-level objective value, and percentage of full LM steps.
		}%
		\label{fig:perf_prof_Ex1}
\end{figure}

For this experiment, we took the maximum number of iteration to be $10^5$ and chose the following parameters for the algorithms: 
$q := 0.8$, 
$\tau_{\textup{abs}} := 10^{-6}$, 
$\tau_{\textup{abs}}^\textup{stat} := 10^{-8}$, 
$\beta := 0.5$, 
$\sigma := 0.5$, 
$\gamma_1 := \gamma_2:=0.5$, 
$\rho_{1} := 10^{-2}$, 
$\rho_{2} := 10^{-12}$, and
$\rho := 10^{-2}$. 
We challenged our algorithms with 121 starting points constructed in the following way. 
The pair $(x,y)$ is chosen from $\{0,1,\ldots,10\}\times\{-5,-4,\ldots,5\}$. 
For the multipliers, we always chose $\mu := \nu := \hat \nu := 1$ and $\lambda := 1$. 
In the case of setting \textbf{Var2}, we made use of  $\zeta:=1$.
The results of the experiments are illustrated in \cref{fig:perf_prof_Ex1} and \cref{tab:data_Ex1}. 

\begin{table}[ht]
    \centering
    \hrule
    \begin{tabular}{c p{1.5cm}p{1.5cm}p{1.5cm}p{1.5cm}p{1.5cm}p{1.5cm}}
     & \textbf{mixLM} & \textbf{FBLM}  & \textbf{mixLM}  & \textbf{FBLM}  & \textbf{mixLM}  & \textbf{FBLM}  \\
     & \textbf{Para}  & \textbf{Para}  & \textbf{Var1}  & \textbf{Var1} & \textbf{Var2}  & \textbf{Var2} \\ 
     \midrule
        Aver.\ \# outer iter.\ & 1476.47 & 1626.61 & 4302.02 & 4301.25 & 1476.47 & 1626.61 \\ 
        Aver.\ \# full  LM steps & 4.3970 & 3.7273 & 3.9256 & 3.1157 & 4.3970 & 3.7273 \\
        Aver.\ time & 0.4748 & 0.4839 & 0.9691 & 1.0176 & 0.4702 & 0.4840 \\
        \# \textbf{Term$=$0} & 13 & 16 & 52 & 52 & 13 & 16 \\
        \# \textbf{Term$=$1} & 63 & 56 & 61 & 49 & 63 & 56 \\
        \# \textbf{Term$=$2} & 45 & 49 & 8 & 20 & 45 & 49 \\
        \# optimal solution	& 74 & 71 & 69 & 68 & 74 & 71 \\  
        \bottomrule
    \end{tabular}
    \caption{%
    	Averaged performance indices for Experiment~1.
    	}%
    \label{tab:data_Ex1}
\end{table}

\begin{figure}[ht]
	\begin{center}
		\includegraphics[height=6cm]{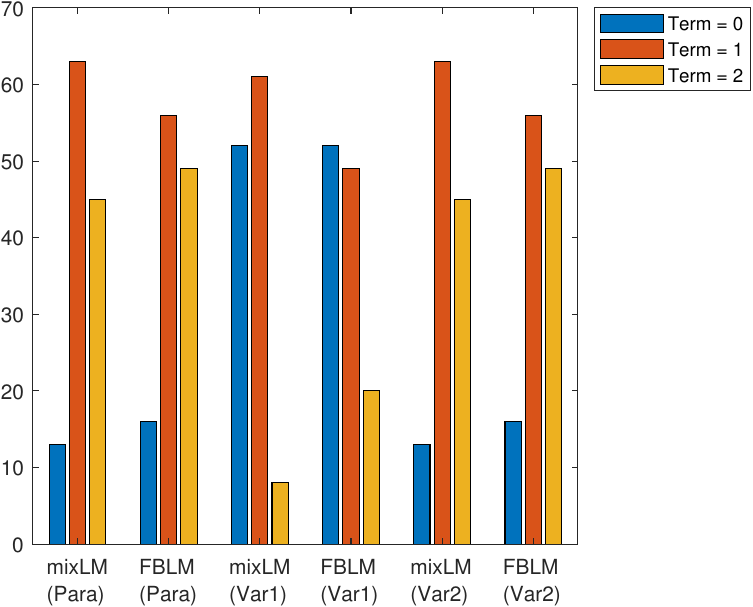} 
	\end{center}
	\caption{Reason for termination for Experiment~1.}\label{fig:termination_Ex1}
\end{figure}

We immediately see that the approaches \textbf{Para} and \textbf{Var2} perform almost equally
good. These methods terminate due to a small residual or stationarity of the squared FB residual
in 90\% of all runs, and the actual optimal solution is found in approximately 60\% of all runs
with slight advantages for \textbf{mixLM}. Let us now report on \textbf{Var1}. The optimal
solution is found in about 57\% of all runs. However, in this setting, both algorithms much more often do not terminate
due to stationarity of the squared FB residual but are aborted since the maximum number
of iterations is reached - the latter happens in 43\% of the runs. Often, both algorithms seem to be close (but not too close)
to a stationary point of the squared FB residual in this case, but this is not detected by our termination
criterion. For a better overview, \cref{fig:termination_Ex1} illustrates the termination behavior of all six methods.
Coming back to an overall comparison, \textbf{mixLM} in setting \textbf{Para} converges to the global minimizer of the
problem for starting points $(x,y)$ chosen from
\[
	(\{0,\ldots,10\}\times\{0,\ldots,5\})\cup(\{5,\ldots,10\}\times\{-2\})\cup(\{9,10\}\times\{-3\}).
\]
In most of the associated runs, a total number of $6$ to $10$ iterations is performed out of which almost all are full
LM steps, i.e., we observe local fast convergence in our experiments. Related phenomena can be observed for the remaining
five approaches. Particularly, despite the absence of any theoretical guarantees, local fast convergence is present in
the variable settings described in \cref{sec:variable_approach}.
	More precisely, we observed that whenever our algorithms approach the global
	minimizer, then this happens in a few number of full LM steps, and the convergence is,
	indeed, quadratic. 
	The large average number of iterations documented in \cref{tab:data_Ex1} results from the
	fact that whenever the algorithms do not approach the global minimizer, they tend to
	get stuck in stationary points of the squared residual which are approached by
	slow gradient steps. Typically, the algorithms are aborted in this situation since
	the maximum number of iterations is reached.
The performance profiles in \cref{fig:perf_prof_Ex1} indicate that settings \textbf{Para} and \textbf{Var2} are superior
to \textbf{Var1} when the total number of iterations and execution time are investigated.
Recalling that \textbf{Var1} does not stop until the maximum number of iterations is reached in comparatively many runs,
its (averaged) running time is twice as large as for \textbf{Para} and \textbf{Var2}.
We note that \textbf{mixLM} (even in setting \textbf{Var1}) performs more full LM steps than \textbf{FBLM}, see \cref{tab:data_Ex1} as well.
For more than 80\% of all starting points, \textbf{mixLM} in setting \textbf{Var1} carries out the highest number of full
LM steps among the six methods which we are comparing here. However, there are some critical instances where, in the setting \textbf{Var1},
comparatively many gradient steps are done, which explains the numbers in \cref{tab:data_Ex1}.
All this is in line with the above comments, and also extends to the computed upper-level objective value, although one has
to be careful as all six methods compute the global minimizer almost equally often. However, \textbf{Var1} ends up in
points with comparatively high objective value in much more runs than the other two settings.
Taking a look at the averaged numbers in \cref{tab:data_Ex1}, we observe that \textbf{mixLM} performs slightly quicker than \textbf{FBLM}
for all three settings. 
Furthermore, \textbf{mixLM} finds the global minimizer more often than \textbf{FBLM}. 
The highest number of total iterations can be observed for \textbf{FBLM} in setting \textbf{Var1}.

\paragraph{Experiment~2}
We consider the bilevel optimization problem
\begin{equation*}
	\min\limits_{x,y}\quad \frac{1}{2}\norm{y - y_\textup{o}}^2\quad\text{s.t.} \quad 
		x \geq 0,\,\texttt{e}^\top x \geq \texttt{e}^\top b^\textup{dem},\,y \in S(x)
\end{equation*}
where $S\colon \mathbb{R}^n \rightrightarrows \mathbb{R}^{n \times \ell}$ 
is the solution mapping of the parametric transportation problem
\begin{equation}\label{eq:parametric_transportation}\tag{TR$(x)$}
	\begin{aligned}
		\min\limits_{y}\quad \sum\limits^{n}_{i=1}\sum\limits_{j=1}^{\ell}c_{ij}y_{ij}
		\quad \text{s.t.} \quad
		\sum\limits_{j=1}^{\ell}y_{ij} &\leq x_{i} & &(i=1,\dots,n),&\\
		\sum\limits_{i=1}^{n}y_{ij} &\geq b_{j}^\textup{dem}&  &(j =1, \dots, \ell),&\\
		y &\geq 0. &&&
	\end{aligned}
\end{equation}
Here, $\ell \in \N$ is a positive integer, $b^\textup{dem}\in \N^{\ell}$ is an integer vector modeling the minimum demand of the $\ell$ consumers, 
and $c \in [0,1]^{n \times \ell}$ is a cost matrix. In \eqref{eq:parametric_transportation}, 
the parameter $x \in \mathbb{R}^n$ represents the offer provided at $n$ warehouses which is unknown 
and shall be reconstructed from a given (noisy) transportation plan $y_\textup{o} \in \mathbb{R}^{n \times \ell}$. 
For our experiments, we chose $n:=5$, $\ell:=7$, and relied on the data matrices given in \cite[Appendix]{Mehlitz2022}.
As mentioned in \cite{Mehlitz2022}, the best known solution of this bilevel optimization problem comes along with
an upper-level objective value of $5.07\cdot 10^{-4}$.

\begin{figure}[ht]
	\begin{center}
		\includegraphics[height=5.0cm]{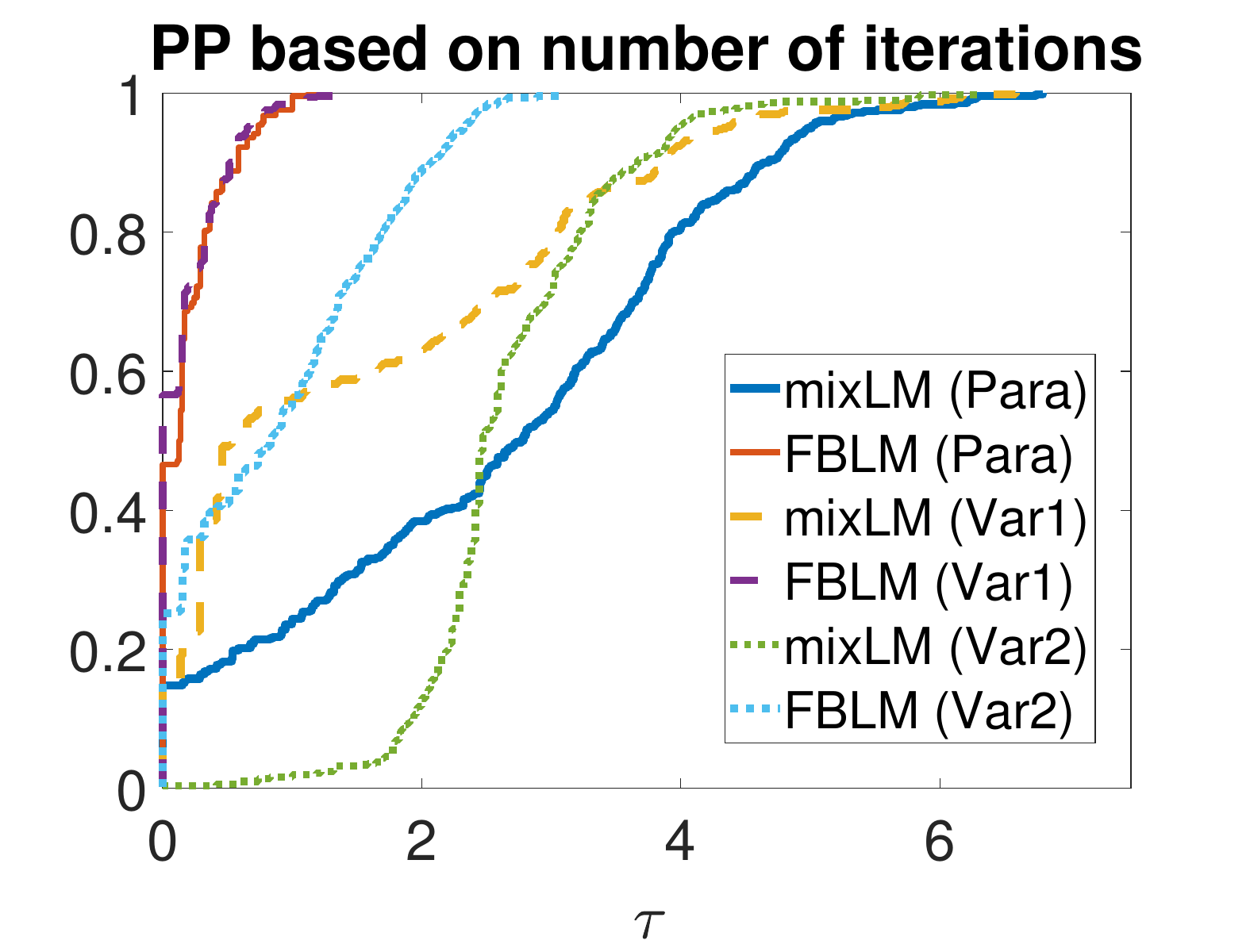} 
		\includegraphics[height=5.0cm]{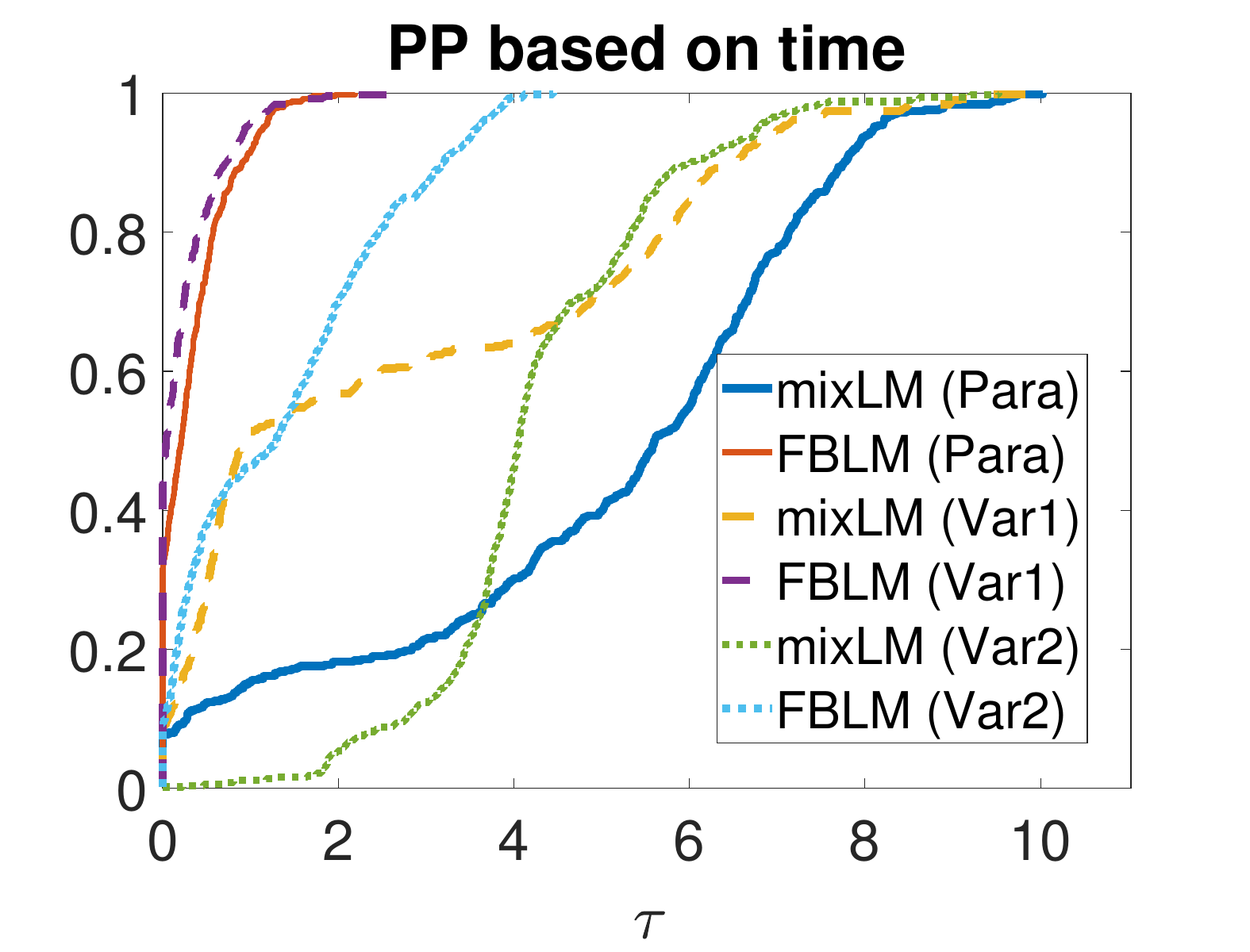}\\ 
		\includegraphics[height=5.0cm]{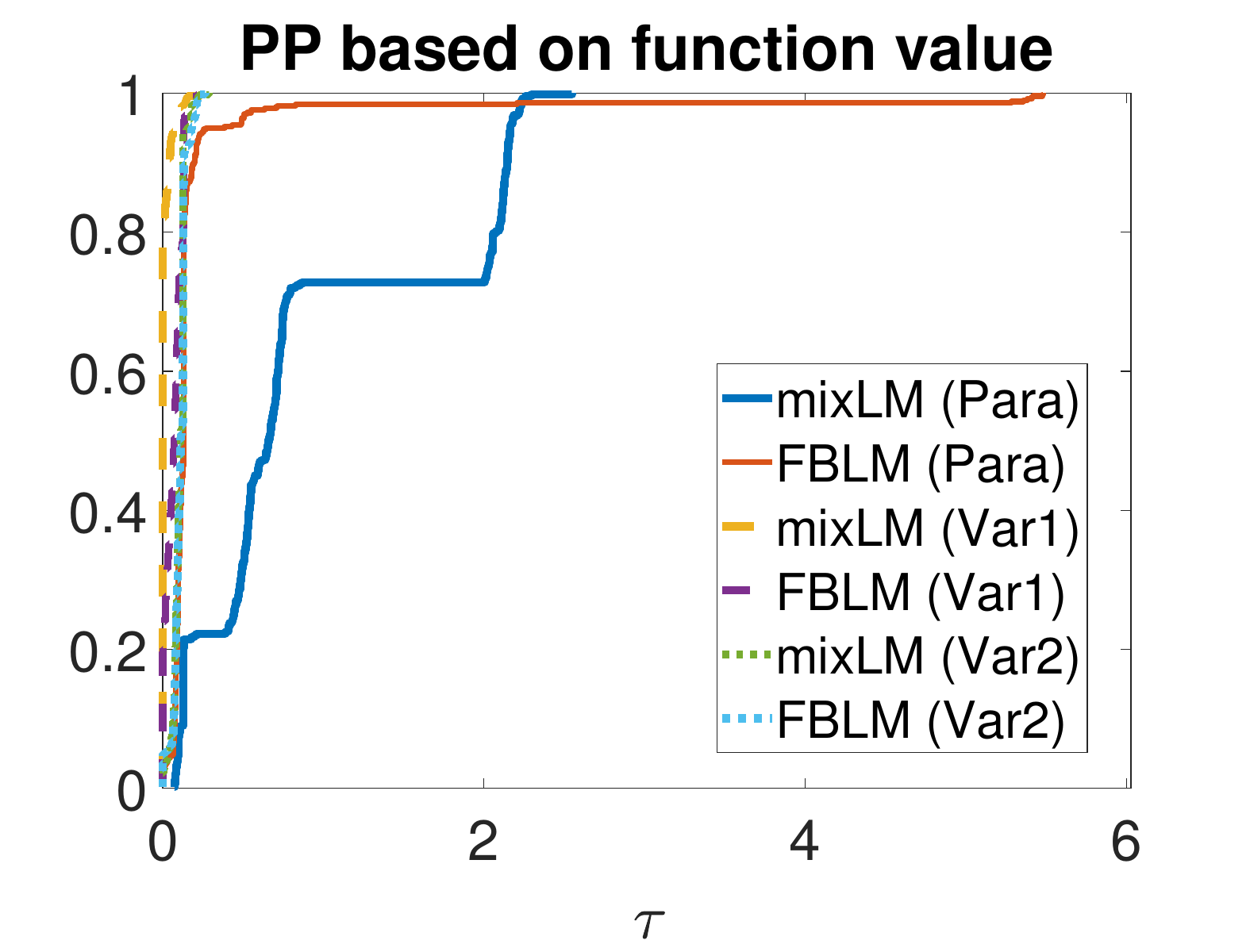} 
		\includegraphics[height=5.0cm]{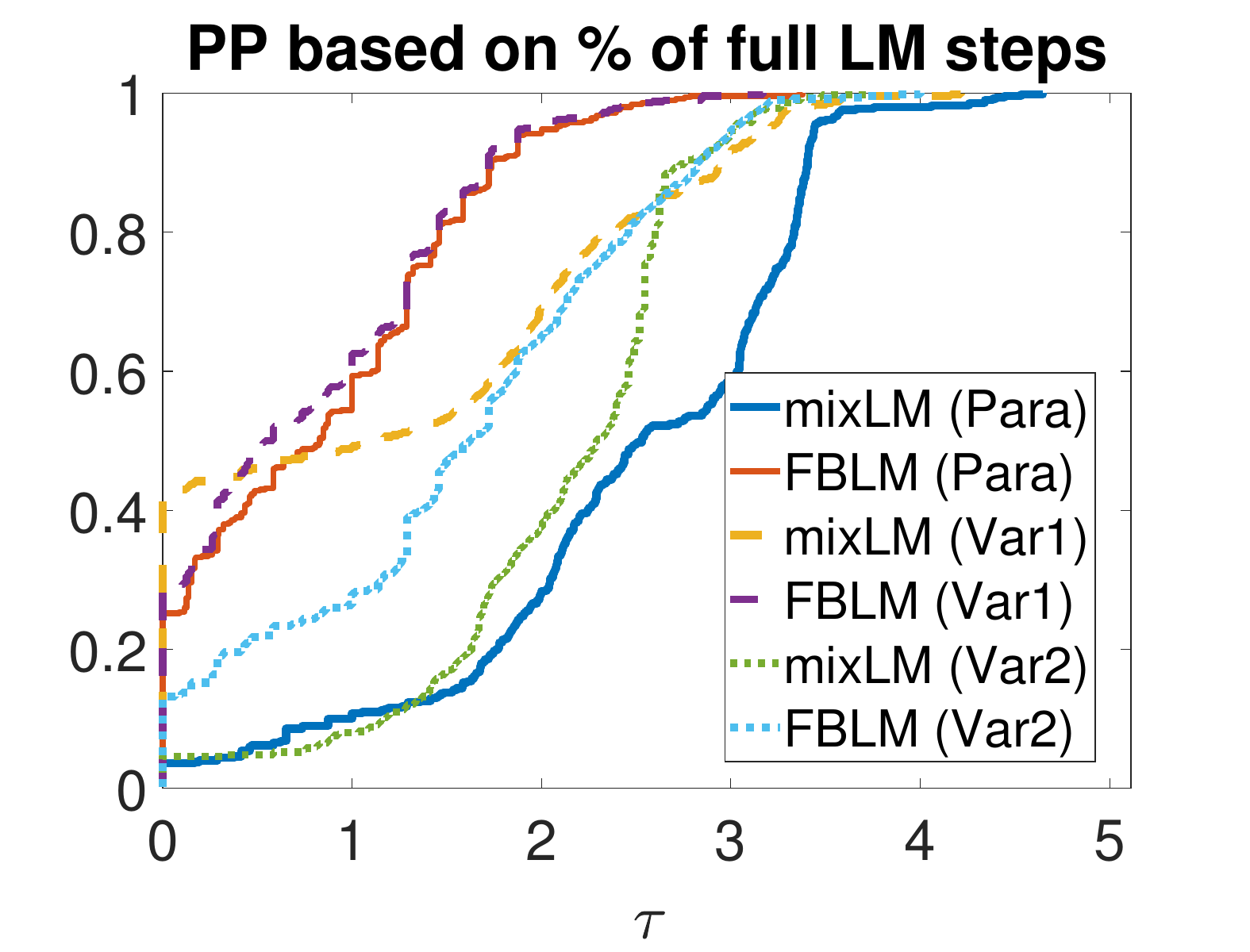} 
	\end{center}
	\caption{%
		Performance profiles for Experiment~2.
		From top left to bottom right: 
    	number of total iterations, execution time, upper-level objective value, and percentage of full LM steps.
		}%
		\label{fig:perf_prof_Ex2}
\end{figure}

All six methods where tested based on a collection of 500 random starting points which were created in the following way. 
For the construction of the pair $(x,y)$, the components of $x$ are chosen randomly from the interval $[1,10]$ 
while the components of $y$ are chosen randomly from $[0,10]$, based on a uniform distribution, respectively.
The associated multiplier vectors as well as $\zeta$ in setting \textbf{Var2} are defined as in our first experiment.
For our algorithms, we chose a maximum number of $10^4$ iterations,
and we made use of the following values for all appearing parameters:
$q := 0.9$, 
$\tau_\textup{abs} := 10^{-4}$, 
$\tau_\textup{abs}^\textup{stat} := 10^{-3}$, 
$\beta := 0.9$, 
$\sigma := 0.4$, 
$\gamma_1 := 10^{-4}$, 
$\gamma_2 := 0.05$, 
and
$\rho_{1} := \rho_2 := \rho := 10^{-4}$.
The resulting performance profiles and averaged performance indices are documented 
in \cref{fig:perf_prof_Ex2} and \cref{tab:data_Ex2}.

\begin{table}[ht]
	\centering
	\hrule
	\begin{tabular}{c p{1.5cm}p{1.5cm}p{1.5cm}p{1.5cm}p{1.5cm}p{1.5cm}}
		 & \textbf{mixLM} & \textbf{FBLM}  & \textbf{mixLM}  & \textbf{FBLM}  & \textbf{mixLM}  & \textbf{FBLM}  \\
     	& \textbf{Para}  & \textbf{Para}  & \textbf{Var1}  & \textbf{Var1} & \textbf{Var2}  & \textbf{Var2} \\ 
     	\midrule
		Aver.\ \# outer iter.\ & 97.754 & 10.256 & 54.508 & 10.156 & 67.75 & 20.17 \\ 
		Aver.\ \# full  LM steps & 5.902 & 7.79 & 8.684 & 7.782 & 21.064 & 8.832 \\
		Aver.\ time & 2.5777 & 0.0391 & 1.0971 & 0.0371 & 0.9610 & 0.1198 \\
		\# \textbf{Term$=$0} & 1 & 0 & 0 & 0 & 0 & 0 \\
		\# \textbf{Term$=$1} & 32 & 301 & 0 & 142 & 0 & 421 \\
		\# \textbf{Term$=$2} & 467 & 199 & 500 & 358 & 500 & 79 \\
		\# approx.\ optimal & 111 & 397 & 0 & 162  & 411 & 378 \\
		\bottomrule
	\end{tabular}
	\caption{Averaged performance indices for Experiment 2.}
	\label{tab:data_Ex2}
\end{table}

\begin{figure}[ht]
	\begin{center}
		\includegraphics[height=6cm]{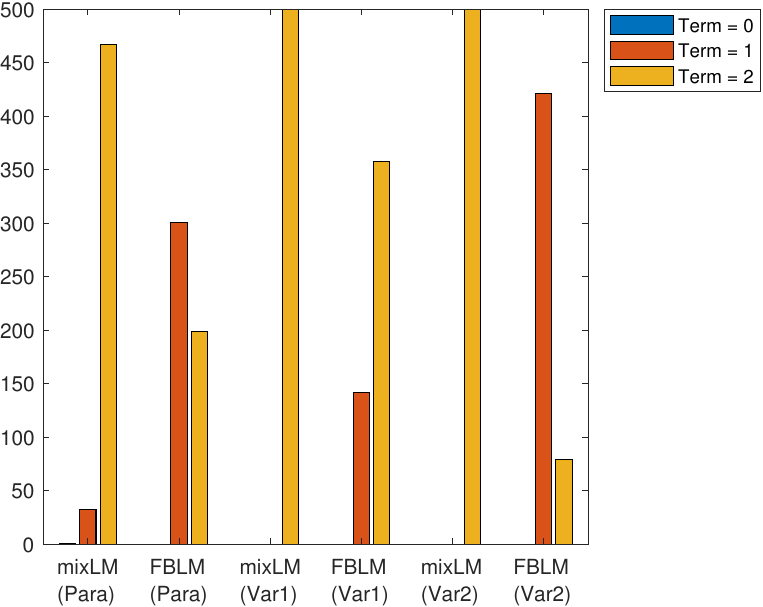} 
	\end{center}
	\caption{Reason for termination for Experiment 2.}\label{fig:termination_Ex2}
\end{figure}

From \cref{fig:perf_prof_Ex2}, it is easy to see that \textbf{FBLM} performs better than \textbf{mixLM} 
considering the total number of iterations, execution time, and the percentage of full LM steps although 
\textbf{FBLM} in setting \textbf{Var2} cannot challenge the same algorithm in settings
\textbf{Para} and \textbf{Var1} for these three criteria, see \cref{tab:data_Ex2} as well.
Our choice of $\tau_\textup{abs}$ and $\tau_\textup{abs}^\textup{stat}$ caused that all six methods terminated
either due to a small residual or stationarity of the squared FB residual in most of the runs. 
Let us point out that \textbf{mixLM} terminated due to stationarity of the squared FB residual
in most of the runs, and the same holds true for \textbf{FBLM} in setting \textbf{Var1}, see \cref{fig:termination_Ex2} 
as well. 

This, however, does not tell the full story. 
In the last row of \cref{tab:data_Ex2}, we document the number of runs where the final iterate produces an upper
level objective value between $5.06\cdot 10^{-4}$ and $5.08\cdot 10^{-4}$, i.e., we count the number of runs
where a point is produced whose associated upper-level function value is close to the best known one 
(below, we will refer to such points as reasonable).
As it turns out, \textbf{mixLM} in setting \textbf{Var1} is not competitive at all in this regard.
Indeed, this method most often terminates in points possessing upper-level objective value around $4.07\cdot 10^{-4}$
and norm of the FB residual around $10^{-5}$. We suspect that this method tends to get stuck in
stationary points of the squared FB residual which are not feasible for the inverse transportation problem.
Hence, the performance profile which monitors the upper-level objective value in \cref{fig:perf_prof_Ex2} is of limited meaningfulness.
The situation is far better for \textbf{FBLM} where reasonable points are computed much more frequently,
although it has to be mentioned that the settings \textbf{Para} (79\% reasonable) 
and \textbf{Var2} (76\% reasonable) again outrun \textbf{Var1} (32\% reasonable) in this regard.
Interestingly, \textbf{mixLM} in setting \textbf{Var2} always terminated due to a small value of the
squared FB residual, but produced reasonable points in more than 82\% of all runs.
An individual fine tuning of the parameters for the settings \textbf{Para}, \textbf{Var1}, and \textbf{Var2}
may lead to more convincing termination reasons, but we abstained from it for a better comparison of all methods
under consideration. 
The performance profiles in \cref{fig:perf_prof_Ex2} show that \textbf{FBLM} in the settings \textbf{Para} and \textbf{Var1} 
are superior to \textbf{Var2}, and when only computed function values are taken into account, then, for the price of
higher computation time, \textbf{mixLM} in setting \textbf{Var2} is acceptable as well.

Let us note that \eqref{eq:KKT_obpp} turns out to be an over-determined mixed linear complementarity system
in the particular setting considered here and, thus, a certain error bound condition is present.
In the light of available literature on smooth LM methods, see e.g.\ \cite{YamashitaFukushima2001}, 
this could be a reason for the observed local fast convergence for a large number of starting points,
although we did not prove such a result in this paper.
Furthermore, it has to be observed that a reformulation via \textbf{Var2} abrogates linearity of the system,
but we obtained the best results in this setting when convergence to reasonable points is the underlying criterion.

\paragraph{Summary}
Let us briefly summarize that our experiments visualized competitiveness of the settings \textbf{Para} and \textbf{Var2},
while setting \textbf{Var1} has turned out to come along with numerical disadvantages in both experiments.
While, in our first experiment, there was no significant difference between \textbf{Para} and \textbf{Var2},
our second (and, by far, more challenging) experiment revealed that \textbf{Var2} also could have some benefits over \textbf{Para}.
Our experiments do not indicate whether it is generally better to apply \textbf{mixLM} or \textbf{FBLM}, but at least we can
suggest that whenever the parameter $\lambda$ is unknown, then it might be reasonable to apply 
\textbf{mixLM} in setting \textbf{Var2} to obtain good solutions.

\section{Concluding remarks}\label{sec:conclusions}

In this paper, we exploited the concept of Newton-differentiability in order to design a nonsmooth
Levenberg--Marquardt method for the numerical solution of over-determined nonsmooth equations.
Our method possesses desirable local convergence properties under reasonable assumptions and
may be applied in situations where the underlying nonsmooth mapping is even discontinuous.
We applied this idea to over-determined mixed nonlinear complementarity systems where a suitable
globalization is possible via gradient steps with respect to the squared norm of the residual
induced by the Fischer--Burmeister function. However, we investigated the method in two different
flavors regarding the computation of the Levenberg--Marquardt direction namely via a residual 
given in terms of the maximum function on the one hand and in terms of the Fischer--Burmeister function on
the other hand. For both methods, global convergence results have been derived, and 
assumptions for local fast convergence have been specified.
Our analysis recovers the impressions from \cite{DeLucaFacchineiKanzow2000}, obtained in a slightly different framework,
that the approach via the maximum residual works under less restrictive assumptions.
Finally, we applied these globalized nonsmooth Levenberg--Marquardt methods in order to solve
bilevel optimization problems numerically. Theoretically, we were in position to verify local fast
convergence of both algorithms under less restrictive assumptions than those ones stated in the literature
which among others assume strict complementarity, see \cite[Section~3]{FliegeTinZemkoho2021}. 
Some numerical experiments were discussed in order
to visualize interesting computational features of this approach.

Our theoretical and numerical investigation of bilevel optimization problems solely focused on the optimistic
approach, but it might be also possible to address (similarly over-determined) stationarity systems related
to pessimistic bilevel optimization via a similar approach.
Next, we would like to point the reader's attention to the fact that the (inherent) inner semicontinuity assumption
on the solution mapping in \cref{lem:CQ_for_stationarity_in_bilevel_programming}, 
see \cref{rem:KKT_obpp} as well, is comparatively strong
and may fail in many practically relevant scenarios.
However, it is well known from the literature that in the presence of so-called inner semicompactness, which is inherent
whenever the solution mapping is locally bounded, slightly more general necessary optimality conditions can
be derived which comprise additional geometric constraints of polyhedral type addressing the multipliers,
see e.g.\ \cite[Theorem~4.9]{DempeZemkoho2011}.
For the numerical solution of this stationarity system, a numerical method would be necessary
which is capable of solving an over-determined system of equations subject to polyhedral constraints 
where either the polyhedral constraint set is not convex or the equations are nonsmooth.
Exemplary, we refer the interested reader to 
\cite{BehlingFischerHaeserRamosSchoenefeld2017,FacchineiFischerHerrich2013,FacchineiFischerHerrich2014} 
for examples of such methods.
In the future, it needs to be checked how these methods apply to the described setting of
bilevel optimization, and how the assumptions for local fast convergence can be verified in this situation.
Yet another way to extend our approach to more general bilevel optimization problems could be to consider the
combined reformulation of bilevel optimization problems which merges the value function and
Karush--Kuhn--Tucker reformulation, see \cite{YeZhu2010}. A potential associated stationarity system can be
reformulated as an overdetermined system of discontinuous nonsmooth equations which, similar as in \cite{HarderMehlitzWachsmuth2021},
still can be solved with the aid of the nonsmooth Levenberg--Marquardt method developed in \cref{sec:LM_local} 
as the latter one is reasonable for functions which are merely calm at their roots, and this property
is likely to hold, see \cite[Lemma~3.3(d)]{HarderMehlitzWachsmuth2021}.
Let us also mention that it would be interesting to study whether an error bound condition is enough to yield
local fast convergence of \cref{alg:GlobalLM,alg:GlobalLM_FB} even in the nonsmooth setting, 
see \cite{YamashitaFukushima2001} for the analysis in the smooth case.
Furthermore, it remains to be seen whether some reasonable conditions can be found which guarantee local fast convergence
of our algorithms when applied to stationarity systems of bilevel optimization problems where the multiplier
associated with the value function constraint is treated as a variable, see \cref{sec:variable_approach}.

\end{document}